\documentclass[11pt]{article}
\usepackage[utf8]{inputenc}
\usepackage[T1]{fontenc}
\usepackage[margin=1in]{geometry}
\usepackage{amsmath,amssymb,amsfonts,amsthm}
\usepackage{booktabs}
\usepackage{enumitem}
\usepackage{hyperref}
\hypersetup{colorlinks=true, linkcolor=blue, urlcolor=blue}
\usepackage{parskip}
\usepackage{mathtools}
\usepackage{graphicx}
\usepackage{float}
\usepackage{caption}

\newtheorem{theorem}{Theorem}[section]
\newtheorem{proposition}[theorem]{Proposition}

\newtheorem{corollary}[theorem]{Corollary}
\newtheorem{lemma}[theorem]{Lemma}
\theoremstyle{definition}
\newtheorem{definition}[theorem]{Definition}
\newtheorem{assumption}[theorem]{Assumption}
\newtheorem{remark}[theorem]{Remark}

\newcommand{\R}{\mathbb{R}}

\newcommand{\Prob}{\mathbb{P}}
\newcommand{\Var}{\mathrm{Var}}
\newcommand{\Cov}{\mathrm{Cov}}
\newcommand{\dto}{\overset{\mathcal{D}}{\longrightarrow}}
\newcommand{\pto}{\overset{\mathbb{P}}{\longrightarrow}}
\DeclareMathOperator*{\argmax}{arg\,max}

\title{\textbf{Asymptotic Divergence-Based Monitoring for Data Drift in
Production AI Systems: A Discrete Case Approach}}
\author{Amadou Diadie Ba \\
LERSTAD, Gaston Berger University, Saint-Louis, Senegal \\
\texttt{ba.amadou-diadie@ugb.edu.sn}}
\date{September 2026}

\begin{document}
\maketitle

\begin{abstract}
We develop a complete asymptotic theory for monitoring data drift in
production artificial intelligence (AI) systems whose inputs take
values in a finite set of categories. The framework combines two
complementary inferential tools. Under the no-drift hypothesis $H_0$,
the two-sample likelihood-ratio statistic $G$ of Wilks follows a
universal $\chi^2_{r-1}$ law, independently of the divergence chosen
and of the sample sizes; the calibrated Population Stability Index
(PSI) of Yurdakul and Naranjo (2021) shares the same $\chi^2_{r-1}$
limit and is asymptotically equivalent to the $G$-test under $H_0$.
Under the drift hypothesis $H_1$, the plug-in $\phi$-divergence is
asymptotically normal with a variance computed in closed form for the
one-sample, two-sample and symmetrized configurations, yielding
standardized $Z$-statistics, detection boundaries, power functions
and sample-size planning rules. We supply a per-category
decomposition that localizes the drift on specific categories, a
Benjamini--Hochberg procedure for monitoring many features
simultaneously, and an explicit comparison with the heuristic PSI
cut-offs of industrial practice. Two fully worked simulation examples
- one one-sample with the R\'enyi divergence and multiple testing,
one two-sample with the Kullback--Leibler (KL) divergence and
sample-size planning - illustrate every notion of the paper on
realistic production data.
\end{abstract}

\vspace{1em}

\vspace{0.5em}

\section{Introduction}

Production machine-learning systems operate on data whose
distribution drifts over time: user populations change, sensors are
recalibrated, upstream pipelines are modified, adversarial behaviour
evolves. Detecting such drift early, and quantifying its magnitude in
a way that supports an operational decision, is a central problem of
AI monitoring. The dominant tool in industrial practice is the
\emph{Population Stability Index} (PSI), used with the heuristic
thresholds $0.10$ and $0.25$. These thresholds have no asymptotic
justification: their Type~I and Type~II errors depend on the number of
categories, on the sample sizes, and on the underlying distributions.
The statistical properties of the PSI have only recently been
established by Yurdakul and Naranjo (2021), who show that, under
$H_0$, $\frac{nm}{n+m}\,\mathrm{PSI}$ follows a $\chi^2_{r-1}$ law
and derive calibrated benchmarks that replace the heuristic cut-offs.

This paper builds on Yurdakul and Naranjo's calibration and extends
it in five directions. First, it covers general $\phi$-divergences
beyond the Kullback--Leibler (KL) divergence: the exact
$\chi^2_{r-1}$ calibration holds for KL, the Jensen--Shannon (JS)
divergence, Tsallis and R\'enyi, while the squared $L_2$ divergence
$L_2^2$ requires a weighted-$\chi^2$ calibration. Second, it provides
the asymptotic distribution of the divergence estimator under $H_1$,
yielding $Z$-statistics, detection boundaries and closed-form power.
Third, it derives the refined symmetrized variance with the
covariance correction needed in the drift regime
$\mathbf p\approx\mathbf q$. Fourth, it adds a per-category
decomposition that localizes the drift on specific categories. Fifth,
it supplies a Benjamini--Hochberg procedure for monitoring many
features simultaneously and a sample-size planning rule. Yurdakul and
Naranjo's calibration is the $H_0$ backbone of the framework; the
present paper is the $H_1$ superstructure and the multi-feature
extension.

The paper is organized as follows. Section~\ref{sec:theory} recalls
the asymptotic theory of $\phi$-divergence estimators.
Section~\ref{sec:detection} introduces the cost of drift, the
directional interpretation, and the per-category decomposition.
Section~\ref{sec:H0} treats the no-drift regime: the hypothesis
formulation, Wilks' theorem, the universal $\chi^2_{r-1}$ calibration,
the one-sample counterpart, and the PSI calibration of Yurdakul and
Naranjo with its equivalence to the $G$-test.
Section~\ref{sec:H1} treats the drift regime: the standardized
$Z$-statistics in the one-sample, two-sample and symmetrized
configurations, their detection boundaries, and their asymptotic
power. Section~\ref{sec:power} develops the power function and the
optimal choice of divergence. Section~\ref{sec:fdr} presents the
Benjamini--Hochberg multiple-testing control.
Section~\ref{sec:simulation-mc} reports a Monte-Carlo simulation study
that validates the asymptotic predictions at finite sample sizes and
quantifies two operational pitfalls.
Section~\ref{sec:simulation} contains the two fully worked simulation
examples on production data. Section~\ref{sec:limitations} discusses
the scope, limitations, and future work.
Section~\ref{sec:guidelines} distils the results into an eight-step
operational checklist. Section~\ref{sec:conclusion} concludes, and the
appendix summarizes the notation and lists the acronyms.

\section{Asymptotic Theory of $\phi$-Divergence Estimators}
\label{sec:theory}

The results of this section are recalled from \cite{BaLo2019}, with
the proofs adapted to the notation of this paper.

\subsection{Notation and definitions}
Let $\mathcal{X} = \{c_1,\ldots,c_r\}$ be a finite set of $r \geq 2$
categories, and let $\mathcal{S} = \{1,\ldots,r\}$ denote the index set.
Let $\mathbf{p} = (p_1,\ldots,p_r)$ and
$\mathbf{q} = (q_1,\ldots,q_r)$ be two probability distributions with
$p_j > 0$ and $q_j > 0$. We observe two independent and identically
distributed (i.i.d.) samples
$X_1,\ldots,X_n \sim \mathbf{p}$ and $Y_1,\ldots,Y_m \sim \mathbf{q}$,
with empirical distributions $\hat{\textbf{p}}_n=\left(\hat{p}_{n,j}\right)_{j\in \mathcal{S}}$ and $\hat{\textbf{q}}_m =\left(\hat{q}_{m,j}\right)_{j\in \mathcal{S}}$ with
\begin{equation}\label{eq:emp-pmf}
\hat{p}_{n,j} = \frac{1}{n}\sum_{i=1}^n \mathbf{1}_{c_j}(X_i),
\qquad
\hat{q}_{m,j} = \frac{1}{m}\sum_{i=1}^m \mathbf{1}_{c_j}(Y_i),
\end{equation}
where $\mathbf{1}_{c_j}(X_i)=1$ if $X_i=c_j$ and $0$ otherwise. The
sum $\sum_{i=1}^n \mathbf{1}_{c_j}(X_i)$ is therefore exactly the
number of observations $X_1,\dots,X_n$ that belong to category $c_j$.
Denoting this count by
\begin{equation}\label{eq:Oj}
O_j := \sum_{i=1}^n \mathbf{1}_{c_j}(X_i),
\qquad j\in\mathcal S,
\end{equation}
we may write $\hat p_{n,j}=O_j/n$: the empirical pmf is simply the
vector of relative frequencies. Note that $\sum_{j\in\mathcal S}O_j=n$,
since each observation belongs to exactly one category.

The formula $\hat p_{n,j}=\frac{1}{n}\sum_{i=1}^n\mathbf{1}_{c_j}(X_i)$
is the one used in \cite{BaLo2019}; the count notation $O_j$ is
introduced here to make the ``number of observations'' reading
explicit, and it will be specialised to $O^{\mathrm{ref}}_j$ and
$O^{\mathrm{prod}}_j$ in Section~\ref{sec:position}.

The same reading applies to the second sample: with
$O'_j := \sum_{i=1}^m \mathbf{1}_{c_j}(Y_i)$, one has
$\hat q_{m,j}=O'_j/m$ and $\sum_{j\in\mathcal S}O'_j=m$.

\paragraph{Gaussian notation.}
Throughout the paper, $\mathcal{N}(\mu,\sigma^2)$ denotes the
univariate normal distribution with mean $\mu$ and variance
$\sigma^2$ and the symbols $\Phi$ and $\Phi^{-1}$ denote the
cumulative distribution function (cdf) and the quantile function of
the standard normal distribution $\mathcal{N}(0,1)$:
\[
\Phi(x) := \int_{-\infty}^{x} \frac{1}{\sqrt{2\pi}}\, e^{-t^2/2}\, dt,
\qquad
z_p := \Phi^{-1}(p), \quad p\in(0,1).
\]
In particular $z_{1-\alpha/2}$ is the $(1-\alpha/2)$-quantile of
$\mathcal{N}(0,1)$, so that $\Prob(Z>z_{1-\alpha/2})=\alpha/2$ for
$Z\sim\mathcal{N}(0,1)$, and $z_{\beta^\star}=\Phi^{-1}(\beta^\star)$
is the $\beta^\star$-quantile used in the sample-size formula of
Proposition~\ref{prop:sample-size}. The symmetry
$\Phi(-x)=1-\Phi(x)$ is used repeatedly in the power computations of
Sections~\ref{sec:H1} and~\ref{sec:power}.

\begin{definition}[$\phi$-divergence]\label{def:phi-div}
For $\phi:(0,1)^2 \to \R$ twice continuously differentiable, the
$\phi$-divergence between two probability mass functions (pmf's)
$\mathbf{p}$ and $\mathbf{q}$ is
\begin{equation}\label{eq:phi-div}
J(\mathbf{p},\mathbf{q}) = \sum_{j\in \mathcal{S}} \phi(p_j, q_j).
\end{equation}
\end{definition}

We use throughout the notations
\begin{equation}\label{eq:phi-derivatives}
\phi_1^{(1)}(s,t):=\frac{\partial\phi}{\partial s}(s,t),
\quad
\phi_2^{(1)}(s,t):=\frac{\partial\phi}{\partial t}(s,t),
\quad
\phi_1^{(2)}(s,t):=\frac{\partial^2\phi}{\partial s^2}(s,t),
\end{equation}
and $\phi_{1,2}^{(2)}(s,t)=\partial^2\phi/\partial s\,\partial t$.

The two first-order partial derivatives are reported in
Table~\ref{tab:divergences} for the five divergences used in this
paper. They are the building blocks of the gradient components
$a_j=\phi_1^{(1)}(p_j,q_j)$ and $b_j=\phi_2^{(1)}(p_j,q_j)$, and of
their reflected counterparts $c_j=\phi_1^{(1)}(q_j,p_j)$ and
$d_j=\phi_2^{(1)}(q_j,p_j)$. The second-order derivative
$\phi_1^{(2)}$ enters the local-expansion coefficient
$\kappa_\phi(p_j)=p_j\,\phi_1^{(2)}(p_j,p_j)$ of
Section~\ref{sec:calibration}.

\begin{assumption}[Regularity of $\phi$]\label{ass:regularity}
There are constants $c_1,c_2>0$ such that $\phi_1^{(1)}(u,u)=c_1$ and
$\phi_2^{(1)}(u,u)=c_2$ for all $u\in(0,1)$.
\end{assumption}

Assumption~\ref{ass:regularity} ensures that if $\mathbf{p}=\mathbf{q}$
then all variances $V_i$ $(i=1,\dots,4)$ defined below vanish, since
for $i\in\{1,3\}$ we have
\[
V_i=\sum_{j\in\mathcal S} p_j c_1^2
- \Bigl(\sum_{j\in\mathcal S} p_j c_1\Bigr)^2=c_1^2(1-1^2)=0,
\]
and for $i\in\{2,4\}$,
\[
V_i=\sum_{j\in\mathcal S} q_j c_2^2
- \Bigl(\sum_{j\in\mathcal S} q_j c_2\Bigr)^2=c_2^2(1-1^2)=0.
\]

\begin{remark}[Non-degeneracy under $H_1$]\label{rem:nondegenerate}
For each of the five divergences of Table~\ref{tab:divergences}, the
condition $V_1=0$ implies $\mathbf p=\mathbf q$. By
Lemma~\ref{lem:compact-var}, $V_1=0$ iff
$a_J=\phi_1^{(1)}(p_J,q_J)$ is $\mathbf p$-a.s.\ constant, i.e.\
$\phi_1^{(1)}(p_j,q_j)=c$ for all $j$ and some $c$. Substituting the
explicit form of $\phi_1^{(1)}$ gives:
\begin{itemize}[leftmargin=*]
\item KL: $1+\log(p_j/q_j)=c$ for all $j$, so $p_j=Kq_j$ with
$K=e^{c-1}$; summing gives $K=1$, hence $\mathbf p=\mathbf q$.
\item JS: $\tfrac12\log(2p_j/(p_j+q_j))=c$ for all $j$ gives
$2p_j/(p_j+q_j)=K$ with $K=e^{2c}$, hence $p_j=Kq_j/(2-K)$; summing
gives $K/(2-K)=1$, so $K=1$ and $p_j=q_j$.
\item Tsallis-$\alpha$: $(\alpha(p_j/q_j)^{\alpha-1}-1)/(\alpha-1)=c$
for all $j$, so $(p_j/q_j)^{\alpha-1}=K$ with
$K=((\alpha-1)c+1)/\alpha$; summing gives $K^{1/(\alpha-1)}=1$, hence
$\mathbf p=\mathbf q$.
\item $L_2^2$: $2(p_j-q_j)=c$ for all $j$; summing gives $c=0$,
hence $\mathbf p=\mathbf q$.
\item R\'enyi-$\alpha$: not separable, but the variance follows from
the Tsallis variance via the delta method, so
$V_{R,\alpha,1}=0\iff V_{T,\alpha,1}=0$, which reduces to the
Tsallis case.
\end{itemize}
Consequently, under $H_1$ all variances are strictly positive and
the $Z$-statistics of Section~\ref{sec:H1} are non-degenerate.
\end{remark}

\begin{remark}[Assumption~\ref{ass:regularity} is not vacuous]
\label{rem:ass-not-vacuous}
The restriction is satisfied by all standard $\phi$-divergences used
in practice, including the five of Table~\ref{tab:divergences} and
the entire class of $f$-divergences
$D_f(\mathbf p\|\mathbf q)=\sum_j q_j f(p_j/q_j)$ with $f(1)=0$, for
which $\phi_1^{(1)}(u,u)=f'(1)$ is constant.
\end{remark}

\begin{table}[H]
\centering
\footnotesize
\renewcommand{\arraystretch}{1.8}
\begin{tabular}{@{}l l l l@{}}
\toprule
Divergence & $\phi(x,y)$ & $\phi_1^{(1)}(x,y)$ & $\phi_2^{(1)}(x,y)$ \\
\midrule
KL & $x\log(x/y)$ & $1+\log(x/y)$ & $-x/y$ \\[2pt]
JS & $\tfrac12 x\log\dfrac{2x}{x+y}+\tfrac12 y\log\dfrac{2y}{x+y}$
   & $\tfrac12\log\dfrac{2x}{x+y}$ & $\tfrac12\log\dfrac{2y}{x+y}$ \\[6pt]
R\'enyi-$\alpha$ & \emph{not separable} & - & - \\[2pt]
Tsallis-$\alpha$ & $\dfrac{x^\alpha y^{1-\alpha}-x}{\alpha-1}$
   & $\dfrac{\alpha (x/y)^{\alpha-1}-1}{\alpha-1}$ & $-x^\alpha y^{-\alpha}$ \\[8pt]
$L_2^2$ & $(x-y)^2$ & $2(x-y)$ & $-2(x-y)$ \\
\bottomrule
\end{tabular}
\caption{The five divergences and their first-order partial
derivatives. R\'enyi is not separable; its variance is obtained from
Tsallis by the delta method (Corollary~\ref{cor:renyi}). JS is the
Jensen--Shannon divergence; it is separable and satisfies
Assumption~\ref{ass:regularity}.}
\label{tab:divergences}
\end{table}

\begin{remark}[Relation to the Cressie--Read power divergence family]
\label{rem:cressie-read}
The family of divergences considered in this paper - KL, JS,
R\'enyi, Tsallis and $L_2^2$ - is a subset of the \emph{power
divergence family} introduced by Cressie and Read (1984) (\cite{Cressie1984}) and studied
extensively by Read and Cressie (1988) (\cite{Read1988}) for discrete goodness-of-fit.
For a single sample and under $H_0$, the asymptotic distribution of
$2m D_\phi(\mathbf p,\hat{\mathbf q}_m)$ is, up to a rescaling, the
one already obtained by Cressie and Read. Our contribution relative
to their work is threefold: (i) a two-sample formulation with the
exact $\chi^2_{r-1}$ calibration via Wilks' theorem
(Section~\ref{sec:H0}); (ii) the refined symmetrized variance with
the covariance correction of Proposition~\ref{prop:sym-var}, which is
specific to the $\mathbf p\approx\mathbf q$ regime and is not in
Cressie and Read; and (iii) the operational layer (detection
boundary, power, sample-size planning, per-category decomposition,
false discovery rate, FDR) developed in
Sections~\ref{sec:H1}--\ref{sec:fdr}.
\end{remark}

\begin{remark}[On the Jensen--Shannon divergence]\label{rem:JS}
The Jensen--Shannon (JS) divergence is bounded,
$0\le D_{JS}(\mathbf p,\mathbf q)\le\log 2$ in nats, which makes its
magnitude directly interpretable on a fixed scale. Its local
behaviour is $D_{JS}\approx\tfrac18\chi^2(\mathbf p,\mathbf q)$ when
$\mathbf p\approx\mathbf q$, and the calibrated statistic
$8m\,D_{JS}$ follows $\chi^2_{r-1}$ under $H_0$. JS is therefore the
natural choice when one needs a divergence whose numerical value can
be compared across datasets without rescaling.
\end{remark}

\begin{remark}[On the Population Stability Index]\label{rem:PSI}
The PSI is defined as
\[
\mathrm{PSI}(\mathbf p,\mathbf q)
=\sum_{j\in\mathcal S}(p_j-q_j)\log\frac{p_j}{q_j}
=\mathrm{KL}(\mathbf p\|\mathbf q)+\mathrm{KL}(\mathbf q\|\mathbf p),
\]
i.e.\ it coincides with the Jeffreys divergence, which is the
symmetrized KL divergence $2J^{(s)}$ of Section~\ref{sec:sym-refined}
with $J=D_{KL}$. It is therefore not listed as a separate divergence
in Table~\ref{tab:divergences}: doing so would double-count the same
quantity.

The PSI is typically used with the following three-tier heuristic
interpretation, reproduced verbatim from industrial practice
(Yurdakul and Naranjo, 2021, p.~89):
\begin{quote}
\emph{$\mathrm{PSI}<0.10$ means a ``little change'',
$0.10\le\mathrm{PSI}<0.25$ means a ``moderate change'', and
$0.25\le\mathrm{PSI}$ means a ``significant change, action
required''.}
\end{quote}
These thresholds have no asymptotic justification: their Type~I and
Type~II errors depend on $r$, on the sample sizes, and on the
underlying distributions. The correct calibration, established by
Yurdakul and Naranjo (2021, Theorem~3.3), depends on the
configuration. In the \emph{two-sample} setting - both $\mathbf p$
and $\mathbf q$ estimated, the standard use case of the PSI in
industry - one has, under $H_0$,
\[
\left(\frac{1}{n}+\frac{1}{m}\right)^{-1}\mathrm{PSI}
\ \dto\ \chi^2_{r-1},
\qquad\text{equivalently}\qquad
\frac{nm}{n+m}\,\mathrm{PSI}\ \dto\ \chi^2_{r-1}.
\]
In the \emph{one-sample} setting - reference $\mathbf p$ fixed, only
$\mathbf q$ estimated, e.g.\ $n\gg m$ - this reduces to
\[
m\,\mathrm{PSI}\ \dto\ \chi^2_{r-1},
\]
consistent with $\kappa_{\mathrm{PSI}}=2$ and
Proposition~\ref{prop:H0-one}. The two-sample formula is the one to
use in practice: the one-sample formula applied to the two-sample
setting inflates the Type~I error by a factor of about $6$ when
$n=m$ (for $r=5$ and $\alpha=0.05$, the actual level becomes
$\approx 0.315$ instead of $0.05$). The present paper extends
Yurdakul and Naranjo's calibration to general $\phi$-divergences
(Section~\ref{sec:calibration}), and complements it with the $H_1$
distribution and power function (Section~\ref{sec:H1}) and with the
per-category decomposition (Section~\ref{sec:position}).
\end{remark}

\begin{table}[H]
\centering
\small
\caption{Heuristic PSI bands used in industrial practice (values
quoted from Yurdakul and Naranjo, 2021) and corresponding calibrated
rule under $H_0$. The calibrated rule depends on $r$, on the sample
sizes and on $\alpha$; it is therefore not reducible to a fixed
grid.}
\label{tab:PSI-bands}
\begin{tabular}{@{}lll@{}}
\toprule
PSI value & Heuristic interpretation & Calibrated rule under $H_0$ (two-sample) \\
\midrule
$\mathrm{PSI}<0.10$ & ``little change''
  & $\frac{nm}{n+m}\mathrm{PSI}\le\chi^2_{r-1,1-\alpha}$: do not reject \\
$0.10\le\mathrm{PSI}<0.25$ & ``moderate change''
  & depends on $r,n,m,\alpha$: may or may not reject \\
$0.25\le\mathrm{PSI}$ & ``significant change''
  & $\frac{nm}{n+m}\mathrm{PSI}>\chi^2_{r-1,1-\alpha}$: reject \\
\bottomrule
\end{tabular}
\end{table}

\begin{remark}[Symmetry, binning, and the fixed-support advantage]
\label{rem:PSI-binning}
Yurdakul and Naranjo (2021, p.~4) note that, in the continuous case,
the PSI is only formally symmetric: the bins are determined by the
percentiles of the base population, so swapping base and target
changes the bin boundaries and hence the value of the PSI. Tom (2023)
goes further and argues that binning is a ``noise generating
irreversible quantization'': the number and boundaries of bins
affect the numerical value of the PSI, and collapsing to a single bin
makes it identically zero. In the present framework the categories
$\mathcal{X}=\{c_1,\ldots,c_r\}$ are fixed a priori (intents,
discrete bins, categories), not derived from the base population, so
the symmetry $\mathrm{PSI}(\mathbf p,\mathbf q)=\mathrm{PSI}(\mathbf
q,\mathbf p)$ holds exactly and neither the binning caveat of
Yurdakul and Naranjo nor the quantization critique of Tom applies.
This is one of the operational advantages of the discrete
fixed-support setting.
\end{remark}

\subsection{First-order asymptotic normality}

Define the gradient components
\begin{equation}\label{eq:grad}
a_j := \phi_1^{(1)}(p_j, q_j),
\qquad
b_j := \phi_2^{(1)}(p_j, q_j), \qquad \forall j\in \mathcal{S}
\end{equation}
and the asymptotic variances
\begin{align}
V_1 &= \sum_{j\in \mathcal{S}} p_j(1-p_j) a_j^2 - 2\sum_{i<j} p_ip_j a_i a_j, \label{eq:V1}\\
V_2 &= \sum_{j\in \mathcal{S}} q_j(1-q_j) b_j^2 - 2\sum_{i<j} q_iq_j b_i b_j. \label{eq:V2}
\end{align}

\begin{lemma}[Compact form of the variance]\label{lem:compact-var}
With $a_j := \phi_1^{(1)}(p_j,q_j)$ as in \eqref{eq:grad},
\begin{equation}\label{eq:V1-compact}
V_1 = \sum_{j\in\mathcal S} p_j a_j^2
- \Bigl(\sum_{j\in\mathcal S} p_j a_j\Bigr)^2
= \mathrm{Var}_{\mathbf p}(a_J),
\end{equation}
where $a_J$ is the random variable taking value $a_j$ with
probability $p_j$. The same identity holds analogously for
\begin{eqnarray}
V_2&=& \sum_{j\in\mathcal S} q_j b_j^2
-\Bigl(\sum_{j\in\mathcal S} q_j b_j\Bigr)^2,
\quad b_j = \phi_2^{(1)}(p_j,q_j),\label{eq:V2'}\\
V_3 &=& \sum_{j\in\mathcal S} q_j c_j^2
-\Bigl(\sum_{j\in\mathcal S} q_j c_j\Bigr)^2,
\quad c_j = \phi_1^{(1)}(q_j,p_j),\label{eq:V3}\\
V_4 & =& \sum_{j\in\mathcal S} p_j d_j^2
-\Bigl(\sum_{j\in\mathcal S} p_j d_j\Bigr)^2,
\quad d_j = \phi_2^{(1)}(q_j,p_j).\label{eq:V4}
\end{eqnarray}
\end{lemma}
\begin{proof}
\emph{Algebraic identity.} Expanding the square,
\[
\Bigl(\sum_j p_j a_j\Bigr)^2 = \sum_j p_j^2 a_j^2
+ 2\sum_{i<j} p_i p_j a_i a_j,
\]
so that
\[
\sum_j p_j a_j^2 - \Bigl(\sum_j p_j a_j\Bigr)^2
= \sum_j p_j(1-p_j)a_j^2 - 2\sum_{i<j} p_i p_j a_i a_j,
\]
which is \eqref{eq:V1}. The right-hand side of \eqref{eq:V1-compact}
is the standard identity
$\mathrm{Var}_{\mathbf p}(a_J) = \mathbb E_{\mathbf p}[a_J^2]
- (\mathbb E_{\mathbf p}[a_J])^2$. The extensions to $V_2$, $V_3$,
$V_4$ follow from the substitutions
$(\mathbf p,a_j)\to(\mathbf q,b_j)$,
$(\mathbf p,a_j)\to(\mathbf q,c_j)$,
$(\mathbf q,b_j)\to(\mathbf p,d_j)$.
\end{proof}

\begin{theorem}[One- and two-sample Central Limit Theorems (CLTs)
{\cite[Thm.~2]{BaLo2019}}, extended]
\label{thm:main}
Assume $p_j > 0$ and $q_j > 0$.
\begin{itemize}
\item[(a)] \emph{One-sample, forward direction.} As $m\to\infty$,
$\sqrt{m}(J(\mathbf{p}, \hat{\mathbf{q}}_m) - J(\mathbf{p},\mathbf{q}))
\dto \mathcal{N}(0, V_2)$.
\item[(b)] \emph{One-sample, reverse direction.} As $m\to\infty$,
$\sqrt{m}(J(\hat{\mathbf{q}}_m, \mathbf{p}) - J(\mathbf{q},\mathbf{p}))
\dto \mathcal{N}(0, V_3)$.
\item[(c)] \emph{Two-sample, forward direction.} If
$n/(n+m)\to\lambda\in(0,1)$,
\begin{equation}\label{eq:two-sample-clt}
\left(\frac{nm}{mV_1+nV_2}\right)^{1/2}
\big(J(\hat{\mathbf{p}}_n, \hat{\mathbf{q}}_m) - J(\mathbf{p},\mathbf{q})\big)
\dto \mathcal{N}(0,1).
\end{equation}
\item[(d)] \emph{Two-sample, reverse direction.} Under the same
condition,
\[
\left(\frac{nm}{nV_3+mV_4}\right)^{1/2}
\big(J(\hat{\mathbf{q}}_m, \hat{\mathbf{p}}_n) - J(\mathbf{q},\mathbf{p})\big)
\dto \mathcal{N}(0,1).
\]
\end{itemize}
\end{theorem}

\paragraph{Reference (parts (a) and (c)).}
Parts (a) and (c) are Theorem~2 of \cite{BaLo2019}; we do not
reproduce the proof here.

\paragraph{Proof of the extension (parts (b) and (d)).}
Parts (b) and (d) follow by the same Taylor-expansion argument,
exchanging the roles of $(\mathbf p,\mathbf q)$ and
$(\mathbf q,\mathbf p)$: the gradient components $a_j,b_j$ are
replaced by $c_j=\phi_1^{(1)}(q_j,p_j)$ and
$d_j=\phi_2^{(1)}(q_j,p_j)$, and the variances $V_1,V_2$ by
$V_3,V_4$. In each case, write
$J(\hat{\mathbf p}_n,\hat{\mathbf q}_m) - J(\mathbf p,\mathbf q)
= \sum_j \bigl[\phi(\hat p_{n,j},\hat q_{m,j})-\phi(p_j,q_j)\bigr]$,
expand each summand to second order around $(p_j,q_j)$, and use the
multinomial CLT for $\sqrt n(\hat{\mathbf p}_n-\mathbf p)$ and
$\sqrt m(\hat{\mathbf q}_m-\mathbf q)$ together with the independence
of the two samples.

\begin{remark}[On the sample-size condition]
In the published version of \cite{BaLo2019}, the condition was stated
as $nm/(n+m)\to\gamma\in(0,1)$, which contains a typo: the quantity
$nm/(n+m)$ tends to $+\infty$ as soon as $n,m\to\infty$. The correct
condition is $n/(n+m)\to\lambda\in(0,1)$, equivalent to
$m/n\to\rho\in(0,\infty)$, and it is the one used throughout this
paper.
\end{remark}
\begin{remark}[Zero probabilities and the choice of divergence]
\label{rem:zeros}
The theory of Section~\ref{sec:theory} assumes $p_j>0$ and $q_j>0$
for all $j\in\mathcal S$. In practice, small categories may have
observed frequencies exactly equal to zero in one of the two samples.
For KL (and hence for the PSI, which is the symmetrized KL), the
plug-in divergence is then undefined, and Tom (2023) documents the
unsatisfactory remedies used in industry (setting the offending term
to zero, or merging bins until both frequencies are positive). The
Jensen--Shannon divergence is immune to this problem: its plug-in
version
$\hat D_{JS}=\frac12\sum_j\bigl[\hat p_j\log(2\hat p_j/(\hat p_j+\hat
q_j))+\hat q_j\log(2\hat q_j/(\hat p_j+\hat q_j))\bigr]$
remains well defined when $\hat p_j=0$ or $\hat q_j=0$, because the
multiplier outside the logarithm vanishes. We therefore recommend JS
as the default divergence when small categories are present, and the
calibrated statistic $8m\,D_{JS}\sim\chi^2_{r-1}$ of
Table~\ref{tab:calibration} applies without modification.
\end{remark}

\begin{corollary}[R\'enyi variances]\label{cor:renyi}
The R\'enyi divergence
$D_{R,\alpha}(\mathbf{p},\mathbf{q})
= \frac{1}{\alpha-1}\log\mathcal{S}_\alpha(\mathbf{p},\mathbf{q})$
with $\mathcal{S}_\alpha = \sum_j p_j^\alpha q_j^{1-\alpha}$ is not
separable. Its asymptotic variances are obtained from the Tsallis
variances by the delta method:
\[
V_{R,\alpha,k} = \frac{V_{T,\alpha,k}}{\mathcal{S}_\alpha^2},
\qquad k\in\{1,2,3,4\}.
\]
\end{corollary}
\begin{proof}
Write $D_{R,\alpha}=g(\mathcal S_\alpha)$ with
$g(s)=\frac{1}{\alpha-1}\log s$. Since
$g'(\mathcal S_\alpha)=1/((\alpha-1)\mathcal S_\alpha)$ and
$D_{T,\alpha}=(\mathcal S_\alpha-1)/(\alpha-1)$, the delta method
gives
$\Var(D_{R,\alpha}) = (g'(\mathcal S_\alpha))^2\Var(\mathcal S_\alpha)
= \Var(D_{T,\alpha})/\mathcal S_\alpha^2$, because
$\Var(\mathcal S_\alpha)=(\alpha-1)^2\Var(D_{T,\alpha})$. The same
rescaling applies to the four orientations $k\in\{1,2,3,4\}$.
\end{proof}

\subsection{Refined variance for symmetrized divergences}
\label{sec:sym-refined}

For a non-symmetric divergence, define the symmetrized version
$J^{(s)}(\mathbf{p},\mathbf{q}) = \tfrac{1}{2}[J(\mathbf{p},\mathbf{q})
+ J(\mathbf{q},\mathbf{p})]$. In the special case $J=D_{KL}$, $2J^{(s)}$
coincides with the PSI (Remark~\ref{rem:PSI}).

\begin{proposition}[Refined symmetrized variance]\label{prop:sym-var}
Let $\eta_j := \tfrac{1}{2}(a_j + d_j)$ and
$\xi_j := \tfrac{1}{2}(b_j + c_j)$. Then the asymptotic variances of
$J^{(s)}(\hat{\mathbf{p}}_n, \mathbf{q})$ and
$J^{(s)}(\mathbf{p}, \hat{\mathbf{q}}_m)$ are
\begin{align}
V_{s,1} &= \tfrac{1}{4}\big(V_1 + V_4 + 2\Cov_{\mathbf{p}}(a_J, d_J)\big), \label{eq:Vs1}\\
V_{s,2} &= \tfrac{1}{4}\big(V_2 + V_3 + 2\Cov_{\mathbf{q}}(b_J, c_J)\big), \label{eq:Vs2}
\end{align}
where
\begin{align*}
\Cov_{\mathbf{p}}(a_J, d_J) &= \sum_j p_j a_j d_j
- \Big(\sum_j p_j a_j\Big)\Big(\sum_j p_j d_j\Big),\\
\Cov_{\mathbf{q}}(b_J, c_J) &= \sum_j q_j b_j c_j
- \Big(\sum_j q_j b_j\Big)\Big(\sum_j q_j c_j\Big).
\end{align*}
\end{proposition}
\begin{proof}
By Taylor expansion of $J^{(s)}$ around $(\mathbf{p}, \mathbf{q})$,
$J^{(s)}(\hat{\mathbf{p}}_n, \mathbf{q}) - J^{(s)}(\mathbf{p},\mathbf{q})
= \sum_j \eta_j (\hat{p}_{n,j} - p_j) + o_P(n^{-1/2})$, because
$\partial J^{(s)}/\partial p_j = \tfrac12(a_j+d_j)=\eta_j$. The
variance of the linear term is
$\Var_{\mathbf{p}}(\eta_J) = \tfrac{1}{4}[V_1 + V_4 + 2\Cov_{\mathbf{p}}(a_J, d_J)]$.
The second identity is analogous, with $\xi_j$ in place of $\eta_j$.
\end{proof}

\begin{remark}[The covariance is not negligible in the drift regime]
\label{rem:cov-negligibility}
In the drift regime $\mathbf{p} \approx \mathbf{q}$, $a_J \approx d_J$
for KL and $L_2^2$, so $\Cov_{\mathbf{p}}(a_J, d_J) \approx V_1$ and
$V_{s,1} \approx V_1$, whereas the uncorrected formula
$(V_1+V_4)/4 \approx V_1/2$. Omitting the covariance therefore
underestimates the variance by a factor of about $2$, and the
standard deviation by a factor of $\sqrt 2$, in the regime where the
framework is actually used. Because PSI coincides with $2J^{(s)}$ for
$J=D_{KL}$, this remark applies to PSI as well: the variance of PSI
cannot be computed by averaging the two KL variances without the
covariance term. All symmetrized computations in this paper use
\eqref{eq:Vs1}--\eqref{eq:Vs2}. This refines the symmetrized
variance formula of \cite[Theorem~4]{BaLo2019}, which was derived
under the implicit assumption that the forward and reverse
linearizations are uncorrelated; in the drift regime
$\mathbf p\approx\mathbf q$, this covariance is of the same order as
the variances themselves and cannot be neglected, as the empirical
study of Section~\ref{sec:mc-sym-var} confirms.
\end{remark}

\begin{theorem}[Symmetrized CLTs {\cite[Thm.~4]{BaLo2019}, refined}]
\label{thm:sym}
Under the assumptions of Theorem~\ref{thm:main}, with $V_{s,1}, V_{s,2}$
as in \eqref{eq:Vs1}--\eqref{eq:Vs2}, as $n\to\infty$,
\[
\sqrt{n}\big(J^{(s)}(\hat{\mathbf{p}}_n, \mathbf{q}) - J^{(s)}(\mathbf{p},\mathbf{q})\big)
\dto \mathcal{N}(0, V_{s,1}),
\]
and if $n/(n+m) \to \lambda \in (0,1)$,
\[
\left(\frac{nm}{mV_{s,1}+nV_{s,2}}\right)^{1/2}
\big(J^{(s)}(\hat{\mathbf{p}}_n, \hat{\mathbf{q}}_m) - J^{(s)}(\mathbf{p},\mathbf{q})\big)
\dto \mathcal{N}(0,1).
\]
\end{theorem}
\begin{proof}
Combine the linear expansion of $J^{(s)}$ in
Proposition~\ref{prop:sym-var} with the multinomial CLT for
$\hat{\mathbf p}_n$ and $\hat{\mathbf q}_m$, exactly as in the proof
of Theorem~\ref{thm:main}.
\end{proof}

\paragraph{From the asymptotic theory to the monitoring framework.}
Section~\ref{sec:theory} has established the asymptotic distribution
of the plug-in $\phi$-divergence under $H_0$ and $H_1$, together with
its variance in the one-sample, two-sample and symmetrized
configurations. We now turn to the monitoring problem itself. The next
section introduces the three notions that operationalise the theory:
the \emph{cost of drift}, which quantifies how far the production
distribution has moved from the reference distribution; the
\emph{per-category contribution}, which decomposes this cost into
interpretable shares that localize the drift on specific categories;
and the \emph{drift position}, which turns these shares into an
actionable shortlist of categories for the operational team. The
statistical tests themselves - the $G$-test under $H_0$ and the
$Z$-statistics under $H_1$ - are developed in
Sections~\ref{sec:H0} and~\ref{sec:H1}, respectively.

\section{Cost of Drift and Drift Position}
\label{sec:detection}

\subsection{Data drift and its informational cost}

\begin{definition}[Data drift]\label{def:drift}
Let $\mathbf p_{\mathrm{train}}$ and $\mathbf q_{\mathrm{prod}}$ be two
pmf's on $\mathcal{X} = \{c_1,\dots,c_r\}$, representing respectively
the reference (training) distribution and the production (monitoring)
distribution. \emph{Data drift} is the event
\[
H_1:\quad \mathbf q_{\mathrm{prod}} \neq \mathbf p_{\mathrm{train}},
\]
i.e.\ there exists at least one category $c_j \in \mathcal{X}$ such
that $q_{\mathrm{prod},j} \neq p_{\mathrm{train},j}$. The no-drift
hypothesis is $H_0:\mathbf q_{\mathrm{prod}} = \mathbf p_{\mathrm{train}}$.
\end{definition}

\begin{definition}[Cost of drift]\label{def:cost}
The \emph{informational cost of drift} between a training pmf
$\mathbf{p}_{\mathrm{train}}$ and a production pmf
$\mathbf{q}_{\mathrm{prod}}$ is
\begin{equation}\label{eq:cost-drift}
D_\phi(\mathbf{p}_{\mathrm{train}}, \mathbf{q}_{\mathrm{prod}})
= \sum_{j \in \mathcal{S}} \phi(p_{\mathrm{train},j}, q_{\mathrm{prod},j}).
\end{equation}
\end{definition}

\begin{remark}[Two notations for the same quantity]
\label{rem:J-vs-D}
The quantities $J(\mathbf p,\mathbf q)$ of
Definition~\ref{def:phi-div} and
$D_\phi(\mathbf p_{\mathrm{train}},\mathbf q_{\mathrm{prod}})$ of
Definition~\ref{def:cost} are the same mathematical object: both equal
$\sum_{j\in\mathcal S}\phi(p_j,q_j)$. We keep two symbols because
they serve two distinct roles - $J$ for the abstract
$\phi$-divergence of Section~\ref{sec:theory}, $D_\phi$ for the
drift-monitoring cost of Section~\ref{sec:detection} - but they are
interchangeable, and every statement proved for one applies verbatim
to the other.
\end{remark}

\paragraph{Directional interpretation.}
For non-symmetric divergences, the direction matters: the forward
divergence
$D_\phi(\mathbf{p}_{\mathrm{train}}, \mathbf{q}_{\mathrm{prod}})$ is
reference-sensitive (it detects the disappearance of training modes),
while the reverse
$D_\phi(\mathbf{q}_{\mathrm{prod}}, \mathbf{p}_{\mathrm{train}})$ is
production-sensitive (it detects the emergence of new modes). In this
paper we adopt the production-first orientation
$D_\phi(\mathbf{q}_{\mathrm{prod}}, \mathbf{p}_{\mathrm{train}})$ as
the default drift indicator, for two reasons. First, it matches the
orientation of the likelihood-ratio $G$-test of Section~\ref{sec:H0}:
the likelihood ratio measures how surprising the production data are
under the reference model, which is exactly the production-first
divergence. Second, it answers the operational question that
motivates monitoring: ``has the production population moved away from
what the model was trained on?''. When no directional preference is
available, the symmetrized divergence $J^{(s)}$ may be used; for KL
this is the PSI (Remark~\ref{rem:PSI}), which is therefore the
symmetric divergence of choice in industrial practice.

\subsection{Per-category contribution and drift position}
\label{sec:position}

The $G$-test of Section~\ref{sec:H0} tells us \emph{whether} a drift
exists, but not \emph{where} it occurs. The present subsection
introduces the per-category decomposition of $G$, whose purpose is to
localize the drift on specific categories and to render the
diagnostic directly actionable.

\begin{remark}[Observed and expected counts]
\label{rem:counts}
In the two-sample setting, we specialise the raw counts $O_j$ and
$O'_j$ of equation~\eqref{eq:Oj} to the reference and production
samples,
\[
O^{\mathrm{ref}}_j := O_j, \qquad
O^{\mathrm{prod}}_j := O'_j, \qquad j\in\mathcal S,
\]
so that $O^{\mathrm{ref}}_j$ ($O^{\mathrm{prod}}_j$) is the number
of observations of $c_j$ in the reference (production) sample, with
$\sum_{j\in\mathcal S}O^{\mathrm{ref}}_j=n$ and
$\sum_{j\in\mathcal S}O^{\mathrm{prod}}_j=m$. The \emph{expected
counts under $H_0$} are defined from the pooled empirical frequency
\[
\hat r_j := \frac{O^{\mathrm{ref}}_j+O^{\mathrm{prod}}_j}{n+m},
\qquad j\in\mathcal S,
\]
by
\[
E^{\mathrm{ref}}_j := n\hat r_j, \qquad
E^{\mathrm{prod}}_j := m\hat r_j, \qquad j\in\mathcal S,
\]
so that $\sum_{j\in\mathcal S}E^{\mathrm{ref}}_j=n$,
$\sum_{j\in\mathcal S}E^{\mathrm{prod}}_j=m$, and
$E^{\mathrm{ref}}_j+E^{\mathrm{prod}}_j=(n+m)\hat r_j$.

In the one-sample setting (Section~\ref{sec:H0-one}), there is no
reference sample, so the reference counts are set formally to
$O^{\mathrm{ref}}_j=E^{\mathrm{ref}}_j:=m\,p_j$; this bookkeeping
convention makes the reference log-ratio
$\log(O^{\mathrm{ref}}_j/E^{\mathrm{ref}}_j)$ vanish identically.
The production counts are $O^{\mathrm{prod}}_j=O'_j$, and the
expected production counts are $E^{\mathrm{prod}}_j:=m\,p_j$. The
signed contribution therefore reduces to
$g_j=O^{\mathrm{prod}}_j\log\bigl(O^{\mathrm{prod}}_j/(m\,p_j)\bigr)$,
and $G_1=2m\,D_{KL}(\hat{\mathbf q}_m,\mathbf p)$.
\end{remark}

\begin{definition}[Per-category contribution and drift position]
\label{def:contribution}
In the two-sample setting with observed counts
$(O^{\mathrm{ref}}_j)_{j\in\mathcal S}$ and
$(O^{\mathrm{prod}}_j)_{j\in\mathcal S}$, and with expected counts
$E^{\mathrm{ref}}_j=n\hat r_j$ and $E^{\mathrm{prod}}_j=m\hat r_j$
under $H_0$, the
\emph{signed contribution} of category $c_j$ to the $G$-statistic is
\[
g_j := O^{\mathrm{ref}}_j \log\frac{O^{\mathrm{ref}}_j}{E^{\mathrm{ref}}_j}
+ O^{\mathrm{prod}}_j \log\frac{O^{\mathrm{prod}}_j}{E^{\mathrm{prod}}_j},
\]
so that $G = 2\sum_{j\in\mathcal S} g_j$, where the observed counts
$O^{\mathrm{ref}}_j$ and $O^{\mathrm{prod}}_j$ are the specialisations
$O_j$ and $O'_j$ of equation~\eqref{eq:Oj} to the reference and
production samples (Remark~\ref{rem:counts}).

The \emph{signed share} of
category $c_j$ is
\[
\pi_j := \frac{g_j}{\sum_{k\in\mathcal S} g_k} = \frac{2 g_j}{G},
\qquad \sum_{j\in\mathcal S} \pi_j = 1.
\]
Its sign indicates the \emph{direction} of the change ($\pi_j>0$:
category $c_j$ is over-represented in production; $\pi_j<0$:
under-represented), and its magnitude indicates the fraction of the
total signal carried by that direction. Because $\pi_j$ can be
negative or exceed $1$ in absolute value when different categories
pull in opposite directions, we also define the \emph{absolute share}
\[
\pi_j^{\mathrm{abs}} := \frac{|g_j|}{\sum_{k\in\mathcal S} |g_k|},
\qquad \pi_j^{\mathrm{abs}}\in[0,1],
\qquad \sum_{j\in\mathcal S} \pi_j^{\mathrm{abs}} = 1,
\]
which measures the fraction of the total signal contributed by
category $c_j$, irrespective of direction. The \emph{position of the
drift} is the subset
\[
\mathcal J^\star := \bigl\{j\in\mathcal S : \pi_j^{\mathrm{abs}} > \gamma\bigr\},
\]
where $\gamma\in(0,1)$ is a user-chosen threshold. In the one-sample
setting (with $\mathbf p$ known), the same construction applies with
$O^{\mathrm{ref}}_j$ formally replaced by $E^{\mathrm{ref}}_j=m\,p_j$:
the signed contribution reduces to
$g_j = O^{\mathrm{prod}}_j \log\bigl(O^{\mathrm{prod}}_j/(m\,p_j)\bigr)$
and $G_1 = 2\sum_j g_j$. The shares $\pi_j$ and
$\pi_j^{\mathrm{abs}}$, and the position $\mathcal J^\star$, are
defined identically.
\end{definition}

\begin{remark}[Why the per-category decomposition is needed]
\label{rem:position-justify}
The $G$-test is a \emph{global} test: it detects whether the joint
distribution has changed, but not \emph{where}. Two very different
drift scenarios can produce the same $G$: a small shift spread across
all $r$ categories (diffuse drift), or a large shift concentrated on
two categories (localized drift). The per-category shares $\pi_j$ and
$\pi_j^{\mathrm{abs}}$ distinguish the two and render the diagnostic
directly actionable, since operational teams can only intervene on
specific categories.
\end{remark}

\begin{remark}[On the default threshold $\gamma = 2/r$]
\label{rem:gamma-default}
The threshold $\gamma$ is not determined by the $G$-test itself: the
test has a canonical threshold $\chi^2_{r-1,1-\alpha}$ under $H_0$,
but the notion of a ``category responsible for the drift'' is a
\emph{reporting} decision rather than a \emph{testing} decision, and
there is no null hypothesis attached to a single category. We
therefore recommend a default and leave the user free to tune it. Our
recommended default is $\gamma = 2/r$. Under a hypothetical uniform
spread of the drift across all $r$ categories, each category would
receive an absolute share of exactly $1/r$. The choice $\gamma = 2/r$
therefore requires that a category contribute \emph{at least twice}
the uniform baseline before being declared part of the drift
position. A larger value of $\gamma$ yields a more conservative
shortlist; a smaller value yields a more inclusive shortlist. The
operational cost of a false omission versus the cost of inspecting
too many categories determines the appropriate value of $\gamma$ in a
given deployment.
\end{remark}

\begin{remark}[Reading the two shares]
\label{rem:signed-vs-absolute}
The signed share $\pi_j$ carries the \emph{direction} of the change
and is the object to report when the operational question is ``which
category has increased and which has decreased?''. The absolute share
$\pi_j^{\mathrm{abs}}$ carries the \emph{magnitude of the
contribution} and is the object to report when the question is
``which categories are driving the alarm?''. Both are needed, but
they play asymmetric roles. Because
$\sum_{k\in\mathcal S}|g_k| \ge G/2 = \sum_{k\in\mathcal S}g_k$, one
always has $\pi_j^{\mathrm{abs}} \le |\pi_j|$: the signed share is
normalized by the total signal $G/2$, which can be small when
categories pull in opposite directions, whereas the absolute share is
normalized by $\sum_k|g_k|$, which cannot be depressed by such
cancellations. When the drift is diffuse and involves categories
pulling in opposite directions, the individual $|\pi_j|$ can
therefore be inflated even though each $\pi_j^{\mathrm{abs}}$ remains
modest; in that case the signed share should not be used for
localization, and the drift position $\mathcal J^\star$ should be
computed from $\pi_j^{\mathrm{abs}}$. Conversely, a category with
$\pi_j^{\mathrm{abs}}$ close to $1$ concentrates the drift on a
single category.
\end{remark}

\section{Under $H_0$: The $G$-Test and its Universal Calibration}
\label{sec:H0}

\subsection{Hypothesis formulation}
\label{sec:hypotheses}

Three equivalent formulations of the null hypothesis arise:
\begin{itemize}[leftmargin=*]
\item \textbf{One-sample with fixed reference:}
$H_0 : D_\phi(\mathbf p, \mathbf q_{\mathrm{prod}}) = 0$,
equivalently $\mathbf q_{\mathrm{prod}} = \mathbf p$.
\item \textbf{Two-sample, both estimated:}
$H_0 : D_\phi(\mathbf p, \mathbf q) = 0$, equivalently
$\mathbf p = \mathbf q$.
\item \textbf{Symmetrized:}
$H_0 : D_\phi^{(s)}(\mathbf p, \mathbf q) = 0$, equivalently
$\mathbf p = \mathbf q$.
\end{itemize}
In all three cases, the variances vanish under $H_0$ by
Remark~\ref{rem:nondegenerate}: $V_1(\mathbf p,\mathbf p) = 0$ and
$V_2(\mathbf p,\mathbf p) = 0$. The $Z$-statistics of
Section~\ref{sec:H1}, which have the form
$Z = \sqrt{N/\hat W}\,\hat D_\phi$ with $\hat W$ a plug-in estimator
of $W$, are therefore the indeterminate form $0/0$ under $H_0$ and
admit no non-degenerate normal limit. The valid test under $H_0$ is
instead the likelihood-ratio $G$-test developed in this section,
whose $\chi^2_{r-1}$ calibration is universal.

Under $H_1$, the variances are strictly positive by
Remark~\ref{rem:nondegenerate}, and the $Z$-statistics are
non-degenerate; this regime is treated in Section~\ref{sec:H1}.

\subsection{Two-sample case}
\label{sec:H0-two}

\subsubsection{The likelihood-ratio statistic and Wilks' theorem}

\begin{theorem}[Wilks, 1938]\label{thm:wilks}
Let $\mathbf{p} = \mathbf{q}$ and $n, m \to \infty$ with
$n/(n+m) \to \lambda \in (0,1)$. Then the likelihood-ratio statistic
\begin{equation}\label{eq:G}
G = 2\Big(n \sum_{j} \hat{p}_{n,j} \log\frac{\hat{p}_{n,j}}{\hat{r}_j}
+ m \sum_{j} \hat{q}_{m,j} \log\frac{\hat{q}_{m,j}}{\hat{r}_j}\Big),
\qquad
\hat r_j = \frac{n\hat p_{n,j}+m\hat q_{m,j}}{n+m},
\end{equation}
converges in distribution to a $\chi^2_{r-1}$ random variable:
$G \dto \chi^2_{r-1}$ as $n, m \to \infty$.
\end{theorem}

\paragraph{Reference.}
Theorem~\ref{thm:wilks} is Wilks' classical theorem
\cite{Wilks1938}. The statistic \eqref{eq:G} is the two-sample
multinomial likelihood ratio for $H_0 : \mathbf{p} = \mathbf{q}$
against $H_1 : \mathbf{p} \neq \mathbf{q}$.

\subsubsection{Decision rule under $H_0$}

\begin{proposition}[Chi-squared rule]\label{prop:H0}
Let $\alpha \in (0,1)$ and let $\chi^2_{r-1,1-\alpha}$ denote the
$(1-\alpha)$-quantile of $\chi^2_{r-1}$. Then, under $H_0$,
\begin{equation}\label{eq:reject-H0}
\Prob\big(G > \chi^2_{r-1,1-\alpha}\big) \to \alpha
\qquad \text{as } n, m \to \infty.
\end{equation}
\end{proposition}
\begin{proof}
Immediate from Theorem~\ref{thm:wilks}.
\end{proof}

\paragraph{Alarm trigger under $H_0$.}
The valid level-$\alpha$ test triggers an alarm whenever
$G > \chi^2_{r-1,1-\alpha}$. The $p$-value is
$p = \Prob(\chi^2_{r-1} > G_{\mathrm{obs}})$. The critical value
$\chi^2_{r-1,1-\alpha}$ is universal: it depends only on $r$ and
$\alpha$, and not on $\mathbf{p}$, $n$, $m$, or $\phi$.

\paragraph{Practical thresholds.}
Table~\ref{tab:chi2-quantiles} reports the universal threshold
$\chi^2_{\nu,1-\alpha}$ for the most common values of $r$ and
$\alpha$.

\begin{table}[H]
\centering
\caption{Universal rejection threshold $\chi^2_{\nu,1-\alpha}$ of
the $G$-test.}
\label{tab:chi2-quantiles}
\begin{tabular}{rcccc}
\toprule
$r$ & $\nu=r-1$ & $\alpha=0.10$ & $\alpha=0.05$ & $\alpha=0.01$ \\
\midrule
2  & 1  & 2.706  & 3.841  & 6.635  \\
3  & 2  & 4.605  & 5.991  & 9.210  \\
5  & 4  & 7.779  & 9.488  & 13.277 \\
10 & 9  & 14.684 & 16.919 & 21.666 \\
20 & 19 & 27.204 & 30.144 & 36.191 \\
\bottomrule
\end{tabular}
\end{table}

\begin{remark}[Rejection, non-rejection, and the misuse of accepting
$H_0$]
\label{rem:accept-H0}
The level $\alpha$ is chosen by the user \emph{before} seeing the
data: $\alpha=0.05$ is the default, $\alpha=0.01$ is used in critical
environments, and $\alpha=0.10$ in exploratory monitoring. Once
$\alpha$ is fixed, the rejection rule is deterministic: reject $H_0$
if and only if $G>\chi^2_{r-1,1-\alpha}$, equivalently if and only if
$p<\alpha$. A non-rejection should not be read as evidence that $H_0$ holds:
whether it is informative depends on the power of the test at the
minimal practically relevant effect, discussed in
Section~\ref{sec:power-role}.
Asymptotic validity of the
$\chi^2_{r-1}$ calibration requires expected counts
$E^{\mathrm{ref}}_j$ and $E^{\mathrm{prod}}_j$ of at least $5$
(Cochran's rule); below this threshold, an exact multinomial test or
a Monte-Carlo calibration should be used instead.
\end{remark}

\subsection{One-sample case with fixed reference}
\label{sec:H0-one}

When $\mathbf{p}$ is known and fixed, the natural one-sample
counterpart of the $G$-test is the KL-based statistic.

\begin{proposition}[One-sample chi-squared limit]\label{prop:H0-one}
Assume $\mathbf{p}$ is known and fixed. Then, under $H_0$ and as
$m \to \infty$,
\begin{equation}\label{eq:chi2-one-sample}
2m\,D_{KL}(\hat{\mathbf{q}}_m, \mathbf{p}) \dto \chi^2_{r-1}.
\end{equation}
The alarm rule under $H_0$ is therefore
$2m\,D_{KL}(\hat{\mathbf{q}}_m, \mathbf{p}) > \chi^2_{r-1,1-\alpha}$,
which has asymptotic level $\alpha$.
\end{proposition}

\begin{proof}
Following the strategy of \cite[Eq.~(5.5)]{BaLo2019}, set
\[
\Delta_j := \hat q_{m,j} - p_j, \qquad
b_m := \max_{j\in\mathcal S}|\Delta_j|,
\qquad j\in\mathcal S.
\]
By the strong law of large numbers, $\Delta_j\to 0$ almost surely
for each $j$, hence $b_m\to 0$ almost surely. Apply the
Taylor--Lagrange--Cauchy formula to the function
$\phi(x,y)=x\log(x/y)$ in its first argument around $p_j$, with
$y=p_j$ fixed. Since $\phi(u,u)=0$, $\phi_1^{(1)}(u,u)=1$ and
$\phi_1^{(2)}(u,u)=1/u$ for all $u\in(0,1)$, we obtain
\[
\phi(\hat q_{m,j},p_j)
= \Delta_j
+ \tfrac12\,\Delta_j^2\,\frac{1}{p_j+\theta_j\Delta_j},
\qquad \theta_j\in(0,1).
\]
Summing over $j$ and multiplying by $2m$,
\begin{equation}\label{eq:KL-decomp}
2m\,D_{KL}(\hat{\mathbf q}_m,\mathbf p)
= 2m\sum_{j\in\mathcal S}\Delta_j
+ m\sum_{j\in\mathcal S}\frac{\Delta_j^2}{p_j+\theta_j\Delta_j}.
\end{equation}
The first term on the right-hand side vanishes identically, because
$\sum_{j\in\mathcal S}\Delta_j=\sum_{j\in\mathcal S}(\hat q_{m,j}-p_j)=1-1=0$.
For the second term, write
\[
\frac{1}{p_j+\theta_j\Delta_j}
= \frac{1}{p_j}
- \frac{\theta_j\Delta_j}{p_j\,(p_j+\theta_j\Delta_j)}.
\]
Since $b_m\to 0$ eventually ensures
$p_j+\theta_j\Delta_j\ge p_j/2$ uniformly in $j$, we have
$\bigl|1/(p_j+\theta_j\Delta_j)-1/p_j\bigr|=O(b_m)$ uniformly in $j$.
Together with $m\sum_j\Delta_j^2/p_j=O_P(1)$ (multinomial CLT) and
$p_j\ge\min_k p_k>0$, this gives
\[
m\sum_{j\in\mathcal S}\Delta_j^2
\Bigl(\frac{1}{p_j+\theta_j\Delta_j}-\frac{1}{p_j}\Bigr)
= o_P(1),
\]
so that \eqref{eq:KL-decomp} becomes
\[
2m\,D_{KL}(\hat{\mathbf q}_m,\mathbf p)
= m\sum_{j\in\mathcal S}\frac{\Delta_j^2}{p_j} + o_P(1).
\]
By the multinomial CLT,
$\sqrt m(\hat{\mathbf q}_m-\mathbf p)\dto\mathcal N(0,\Sigma_{\mathbf p})$
with $\Sigma_{\mathbf p}=\mathrm{diag}(\mathbf p)-\mathbf p\mathbf p^\top$,
and Cochran's theorem gives
$m\sum_{j\in\mathcal S}\Delta_j^2/p_j\dto\chi^2_{r-1}$.
Slutsky's lemma then yields
$2m\,D_{KL}(\hat{\mathbf q}_m,\mathbf p)\dto\chi^2_{r-1}$.
\end{proof}

\subsection{PSI calibration and its equivalence with the $G$-test}
\label{sec:PSI-G}

The two-sample $G$-test and the calibrated PSI test are
asymptotically equivalent under $H_0$. We record this fact as a
separate remark, since it unifies the likelihood-ratio viewpoint of
Section~\ref{sec:H0-two} with the divergence viewpoint of
Remark~\ref{rem:PSI}.

\begin{remark}[PSI calibration and the $G$-test]\label{rem:PSI-G}
Under $H_0$, the $G$-statistic of Theorem~\ref{thm:wilks} and the
calibrated PSI statistic $\frac{nm}{n+m}\,\mathrm{PSI}$ are
asymptotically equivalent: both converge to $\chi^2_{r-1}$. Indeed,
a second-order Taylor expansion of
$G=2[nD_{KL}(\hat{\mathbf p}_n\|\hat{\mathbf r})
+mD_{KL}(\hat{\mathbf q}_m\|\hat{\mathbf r})]$ around
$\hat{\mathbf r}$ gives
\[
G\approx\frac{nm}{n+m}\sum_j\frac{(\hat p_{n,j}-\hat q_{m,j})^2}{\hat r_j},
\]
which coincides with $\frac{nm}{n+m}\mathrm{PSI}^{**}$ at leading
order (Yurdakul and Naranjo, 2021, Theorem~3.2). The two tests
therefore share the same level-$\alpha$ rejection region up to
$o_P(1)$, and the $G$-test can be interpreted as the
likelihood-ratio counterpart of the calibrated PSI. They differ under
$H_1$: the $G$-test is locally most powerful against contiguous
alternatives (Neyman--Pearson), whereas the PSI test has a different
power profile, studied numerically by Yurdakul and Naranjo (2021,
Table~4).
\end{remark}

\begin{table}[H]
\centering
\small
\caption{Calibrated two-sample PSI threshold
$\left( \frac{1}{n}+\frac{1}{m}\right)\chi^2_{r-1,1-\alpha}$ for $r=10$, $\alpha=0.05$
(reproduced from Yurdakul and Naranjo, 2021, Table~2). The heuristic
thresholds $0.10$ and $0.25$ coincide with the calibrated threshold
only in narrow ranges of $(n,m)$.}
\label{tab:PSI-benchmarks}
\begin{tabular}{rcccccc}
\toprule
$n\backslash m$ & $100$ & $200$ & $400$ & $600$ & $800$ & $1000$ \\
\midrule
$100$  & $0.338$ & $0.254$ & $0.211$ & $0.197$ & $0.190$ & $0.186$ \\
$200$  & $0.254$ & $0.169$ & $0.127$ & $0.113$ & $0.106$ & $0.102$ \\
$400$  & $0.211$ & $0.127$ & $0.085$ & $0.070$ & $0.063$ & $0.059$ \\
$600$  & $0.197$ & $0.113$ & $0.070$ & $0.056$ & $0.049$ & $0.045$ \\
$800$  & $0.190$ & $0.106$ & $0.063$ & $0.049$ & $0.042$ & $0.038$ \\
$1000$ & $0.186$ & $0.102$ & $0.059$ & $0.045$ & $0.038$ & $0.034$ \\
\bottomrule
\end{tabular}
\end{table}

The calibration of the four other divergences under $H_0$ is deferred
to Section~\ref{sec:calibration}, where it is combined with the power
analysis to guide the choice of $\phi$.

\section{Under $H_1$: The $Z$-Statistics and Their Detection
Boundaries}
\label{sec:H1}

Under $H_1$, $D_\phi(\mathbf{p}, \mathbf{q}) > 0$ and the variances
are strictly positive. Three configurations arise; we derive the
standardized statistic $Z$, its asymptotic distribution, the
detection boundary $\tau_\alpha$, and the power in each case.

\begin{remark}[Interpretation of the detection boundary]
\label{rem:tau-interpret}
In what follows, $\tau_\alpha$ denotes a \emph{power-analysis
boundary} or \emph{minimum detectable effect}: the threshold such
that the power of the $Z$-test to detect a drift amplitude $D_\phi$
exceeds $1-\alpha/2$ whenever $D_\phi > \tau_\alpha$. It is \emph{not}
a valid level-$\alpha$ test under $H_0$, because under $H_0$ the
variance $W$ vanishes and $Z$ is degenerate. The valid level-$\alpha$
test under $H_0$ remains the universal $G$-test of
Section~\ref{sec:H0}.
\end{remark}

\subsection{Generic structure and assumptions}
\label{sec:H1-generic}

\paragraph{The standardized statistic.}
Let $\hat\theta$ be an estimator of a scalar quantity $\theta$ based
on a sample of size $N$, and let $\hat\sigma_{\hat\theta}$ be a
plug-in estimator of the standard error of $\hat\theta$. The
standardized statistic is $t = (\hat\theta-\theta)/\hat\sigma_{\hat\theta}$.
Under a fixed alternative, $t$ is approximately normal with mean
equal to its \emph{signal-to-noise ratio} (SNR) $\mathrm{SNR}=\theta/\sigma_{\hat\theta}$
and variance one.

In this paper the estimator is the plug-in divergence $\hat D_\phi$.
Three variance-related quantities must be distinguished:
\begin{itemize}[leftmargin=*]
\item the asymptotic variance of the \emph{normalized} estimator
$\sqrt N\,\hat D_\phi$, denoted
$W := \lim_{N\to\infty} N\,\Var(\hat D_\phi)$; this quantity does not
depend on $N$, and it enters the CLTs of Theorem~\ref{thm:main};
\item the asymptotic variance of the \emph{unnormalized} estimator
$\hat D_\phi$, which equals $W/N$;
\item the standard error of $\hat D_\phi$, equal to $\sqrt{W/N}$,
whose plug-in estimator is $\sqrt{\hat W/N}$.
\end{itemize}
The standardized statistic is therefore
\[
Z = \frac{\hat D_\phi}{\sqrt{\hat W/N}} = \sqrt{\frac{N}{\hat W}}\,\hat D_\phi.
\]

\begin{definition}[Standardized statistic]\label{def:Z-generic}
In each configuration,
\begin{equation}\label{eq:Z-generic}
Z = \sqrt{N/\hat W}\, \hat D_\phi,
\end{equation}
where $\hat D_\phi$ is the plug-in divergence, $N$ is a
normalization rate, and $\hat W$ is a weakly consistent plug-in
estimator of $W = \Var(\sqrt{N}\hat D_\phi)$, i.e.\ $\hat W \pto W > 0$.
\end{definition}

\begin{proposition}[Generic asymptotic normality]\label{prop:Z-generic}
Under $H_1$ and the assumptions of Theorem~\ref{thm:main}, for each
fixed pair $(n, m)$ with $n, m$ large,
\begin{equation}\label{eq:Z-normal}
Z \approx \mathcal{N}(\mathrm{SNR}_\phi, 1),
\qquad
\mathrm{SNR}_\phi := \sqrt{N/W}\, D_\phi(\mathbf{p}, \mathbf{q}).
\end{equation}
\end{proposition}
\begin{proof}
By the CLT of Theorem~\ref{thm:main},
$\sqrt{N/W}\,(\hat D_\phi - D_\phi) \dto \mathcal{N}(0, 1)$.
Since $\hat W \pto W$, Slutsky's lemma gives
$\sqrt{N/\hat W}\,(\hat D_\phi - D_\phi) \dto \mathcal{N}(0, 1)$.
Decompose
$Z = \sqrt{N/\hat W}\,(\hat D_\phi - D_\phi) + \sqrt{N/\hat W}\, D_\phi$.
The first term tends to $\mathcal{N}(0,1)$; the second equals
$\sqrt{N/W}\,D_\phi = \mathrm{SNR}_\phi$ by Slutsky again.
\end{proof}

\begin{proposition}[Generic detection boundary]\label{prop:tau-generic}
Let $\alpha \in (0,1)$. Under $H_1$, the condition
$|Z| > z_{1-\alpha/2}$ is equivalent, for $\hat D_\phi \geq 0$, to
\begin{equation}\label{eq:tau-generic}
\hat D_\phi > \tau_\alpha,
\qquad
\tau_\alpha := z_{1-\alpha/2}\, \sqrt{\hat W/N}.
\end{equation}
\end{proposition}
\begin{proof}
Since $Z \geq 0$ under $H_1$, $|Z| = Z$; substitute
\eqref{eq:Z-generic} and multiply by $\sqrt{\hat W/N} > 0$.
\end{proof}

\begin{proposition}[Generic power]\label{prop:power-generic}
Under $H_1$ and the assumptions of Theorem~\ref{thm:main}, the power
of the two-sided level-$\alpha$ test $|Z|>z_{1-\alpha/2}$ is
\begin{equation}\label{eq:power-generic}
\beta_\phi \approx
\Phi\bigl(\mathrm{SNR}_\phi - z_{1-\alpha/2}\bigr)
+ \Phi\bigl(-\mathrm{SNR}_\phi - z_{1-\alpha/2}\bigr),
\end{equation}
where $\mathrm{SNR}_\phi = \sqrt{N/W}\,D_\phi(\mathbf p,\mathbf q)$.
\end{proposition}
\begin{proof}
By Proposition~\ref{prop:Z-generic},
$Z \approx \mathcal N(\mathrm{SNR}_\phi,1)$ under $H_1$. The test
rejects when $|Z|>z_{1-\alpha/2}$, hence
\begin{align*}
\beta_\phi
&= \Prob\bigl(|Z|>z_{1-\alpha/2}\bigr)\\
&= \Prob\bigl(Z>z_{1-\alpha/2}\bigr)
 + \Prob\bigl(Z<-z_{1-\alpha/2}\bigr)\\
&= \Phi\bigl(\mathrm{SNR}_\phi - z_{1-\alpha/2}\bigr)
 + \Phi\bigl(-\mathrm{SNR}_\phi - z_{1-\alpha/2}\bigr),
\end{align*}
using the symmetry $\Phi(-x)=1-\Phi(x)$.
\end{proof}

\begin{remark}[Calibration versus ad-hoc thresholds]
\label{rem:ad-hoc}
Unlike the ad-hoc thresholds used in classical drift detection
(Drift Detection Method, DDM; Adaptive Windowing, ADWIN;
Page--Hinkley) and in the heuristic PSI cut-offs of $0.10$ and
$0.25$, the threshold $\tau_\alpha$ is automatically calibrated by
the CLT of Theorem~\ref{thm:main}: it depends only on $n$, $m$, and
the estimated variances. No subjective choice of a magic number such
as $\tau=0.10$ is required. The level $\alpha$ is the only user-set
quantity, and it has a precise asymptotic interpretation as the
Type~I error of the associated test.
\end{remark}

\begin{remark}[The $G$-test under $H_1$ and the role of the $Z$-statistic]
\label{rem:G-vs-Z}
Under $H_1$, the $G$-test of Theorem~\ref{thm:wilks} is also
consistent, and by Wilks' noncentral theorem $G$ follows
approximately a noncentral $\chi^2_{r-1}$ with noncentrality
parameter
$\Lambda = \frac{nm}{n+m}\sum_j (p_j-q_j)^2/p_j$ under contiguous
alternatives. By the Neyman--Pearson lemma, the likelihood ratio $G$
is locally most powerful against contiguous alternatives; in the
diffuse regime $\mathbf p\approx\mathbf q$, no $Z$-statistic of a
$\phi$-divergence can have higher power than $G$. The $Z$-statistics
of this section are therefore \emph{complementary} to $G$, not
substitutes for it: they quantify (i) the effect size in nats, (ii)
the sensitivity of a specific divergence to a specific drift
direction, and (iii) the per-category localization. A practitioner
who only needs to decide ``drift or not'' should use $G$; a
practitioner who needs to \emph{quantify}, \emph{localize}, or
\emph{prioritize} should use the $Z$-statistics in addition to $G$.
\end{remark}

\subsection{Configuration 1: one-sample with fixed reference}
\label{sec:H1-1s}

\paragraph{Setup.}
The reference $\mathbf p$ is known (from a large historical sample);
the production $\mathbf q$ is estimated by the empirical pmf
$\hat{\mathbf q}_m$. In this configuration we adopt the
production-first orientation $D_\phi(\hat{\mathbf q}_m, \mathbf p)$,
which matches the orientation of the $G$-test and answers the
operational question. Its variance is $V_3$ by
Theorem~\ref{thm:main}(b). Set
$\hat D_\phi = D_\phi(\hat{\mathbf q}_m, \mathbf p)$, $N = m$,
$W = V_3$.

\paragraph{Standardized statistic.}
\begin{equation}\label{eq:Z1}
Z_1 = \frac{\sqrt{m}\, D_\phi(\hat{\mathbf q}_m, \mathbf p)}{\sqrt{\hat V_3}},
\qquad
\hat V_3 = \sum_j \hat{q}_{m,j}\, \hat c_j^2
- \Big(\sum_j \hat{q}_{m,j}\, \hat c_j\Big)^2,
\quad \hat c_j := \phi_1^{(1)}(\hat{q}_{m,j}, p_j).
\end{equation}

\paragraph{Asymptotic distribution under $H_1$.}
By the one-sample CLT,
$Z_1 \approx \mathcal{N}\big(\mathrm{SNR}_\phi^{(1)}, 1\big)$ with
\[
\mathrm{SNR}_\phi^{(1)} := \frac{\sqrt{m}\, D_\phi(\mathbf q, \mathbf p)}{\sqrt{V_3}}.
\]

\paragraph{Detection boundary and power.}
\[
\tau_\alpha^{(1)}(m) = z_{1-\alpha/2}\, \sqrt{\hat V_3 / m},
\qquad
\beta^{(1)} \approx \Phi(\mathrm{SNR}^{(1)}-z_{1-\alpha/2})
+ \Phi(-\mathrm{SNR}^{(1)}-z_{1-\alpha/2}),
\]
by Proposition~\ref{prop:tau-generic} and
Proposition~\ref{prop:power-generic} evaluated at
$\mathrm{SNR}_\phi^{(1)}$.

\subsection{Configuration 2: two-sample with both estimated}
\label{sec:H1-2s}

\paragraph{Setup.}
Both $\mathbf{p}$ and $\mathbf{q}$ are estimated by their empirical
pmf's. Set
$\hat D_\phi = D_\phi(\hat{\mathbf{p}}_n, \hat{\mathbf{q}}_m)$,
$N = nm$, $W = mV_1 + nV_2$.

\paragraph{Standardized statistic.}
\begin{equation}\label{eq:Z2}
Z_2 = \left(\frac{nm}{m\hat V_1 + n\hat V_2}\right)^{1/2}
D_\phi(\hat{\mathbf{p}}_n, \hat{\mathbf{q}}_m),
\end{equation}
where
\begin{align}
\hat V_1 &= \sum_j \hat{p}_{n,j}\, \hat a_j^2
- \Big(\sum_j \hat{p}_{n,j}\, \hat a_j\Big)^2,
& \hat a_j &:= \phi_1^{(1)}(\hat{p}_{n,j}, \hat{q}_{m,j}), \\
\hat V_2 &= \sum_j \hat{q}_{m,j}\, \hat b_j^2
- \Big(\sum_j \hat{q}_{m,j}\, \hat b_j\Big)^2,
& \hat b_j &:= \phi_2^{(1)}(\hat{p}_{n,j}, \hat{q}_{m,j}).
\end{align}

\paragraph{Asymptotic distribution, boundary, power.}
\[
Z_2 \approx \mathcal{N}\!\big(\mathrm{SNR}_\phi^{(2)}, 1\big),
\qquad
\mathrm{SNR}_\phi^{(2)} := \frac{\sqrt{nm}\, D_\phi(\mathbf{p}, \mathbf{q})}{\sqrt{mV_1 + nV_2}},
\]
\[
\tau_\alpha^{(2)}(n,m) = z_{1-\alpha/2}\,
\sqrt{\frac{m\hat V_1 + n\hat V_2}{nm}},
\qquad
\beta^{(2)} \approx \Phi(\mathrm{SNR}^{(2)}-z_{1-\alpha/2})
+ \Phi(-\mathrm{SNR}^{(2)}-z_{1-\alpha/2}),
\]
by Proposition~\ref{prop:tau-generic} and
Proposition~\ref{prop:power-generic} evaluated at
$\mathrm{SNR}_\phi^{(2)}$.

\subsection{Configuration 3: symmetrized with refined variance}
\label{sec:H1-sym}

\paragraph{Setup.}
For a non-symmetric divergence, the direction of the drift is a
modelling choice: forward $D_\phi(\mathbf p,\mathbf q)$ is
reference-sensitive, reverse $D_\phi(\mathbf q,\mathbf p)$ is
production-sensitive, and neither is canonical when the user has no
directional preference. The symmetrized divergence
\begin{equation}\label{eq:Js-def}
J^{(s)}(\mathbf p,\mathbf q) = \tfrac12\bigl[J(\mathbf p,\mathbf q)+J(\mathbf q,\mathbf p)\bigr]
\end{equation}
removes the choice. We estimate it by
$\hat D_\phi^{(s)} = J^{(s)}(\hat{\mathbf p}_n,\hat{\mathbf q}_m)$
and set $N=nm$, $W=mV_{s,1}+nV_{s,2}$, where $V_{s,1},V_{s,2}$ are
the refined variances of Proposition~\ref{prop:sym-var}.

\paragraph{Standardized statistic.}
\begin{equation}\label{eq:Zs}
Z_{n,m}^{(s)} = \Bigl(\frac{nm}{m\hat V_{s,1}+n\hat V_{s,2}}\Bigr)^{1/2}
J^{(s)}(\hat{\mathbf p}_n,\hat{\mathbf q}_m),
\end{equation}
where
\begin{equation}\label{eq:Vs-hat}
\hat V_{s,1}=\tfrac14\bigl(\hat V_1+\hat V_4 + 2\widehat{\Cov}_{\mathbf p}(a_J,d_J)\bigr),
\quad
\hat V_{s,2}=\tfrac14\bigl(\hat V_2+\hat V_3 + 2\widehat{\Cov}_{\mathbf q}(b_J,c_J)\bigr),
\end{equation}
and the plug-in covariances are
\begin{align}
\widehat{\Cov}_{\mathbf p}(a_J,d_J) &=\sum_j\hat p_{n,j}\hat a_j\hat d_j
-\Bigl(\sum_j\hat p_{n,j}\hat a_j\Bigr)\Bigl(\sum_j\hat p_{n,j}\hat d_j\Bigr),\\
\widehat{\Cov}_{\mathbf q}(b_J,c_J) &=\sum_j\hat q_{m,j}\hat b_j\hat c_j
-\Bigl(\sum_j\hat q_{m,j}\hat b_j\Bigr)\Bigl(\sum_j\hat q_{m,j}\hat c_j\Bigr),
\end{align}
with $\hat a_j=\phi_1^{(1)}(\hat p_{n,j},\hat q_{m,j})$,
$\hat b_j=\phi_2^{(1)}(\hat p_{n,j},\hat q_{m,j})$,
$\hat c_j=\phi_1^{(1)}(\hat q_{m,j},\hat p_{n,j})$,
$\hat d_j=\phi_2^{(1)}(\hat q_{m,j},\hat p_{n,j})$.

\paragraph{Asymptotic distribution, boundary, power.}
By Theorem~\ref{thm:sym},
\[
Z_{n,m}^{(s)}\approx\mathcal N\bigl(\mathrm{SNR}_\phi^{(s)},1\bigr),
\qquad
\mathrm{SNR}_\phi^{(s)} := \frac{\sqrt{nm}\,J^{(s)}(\mathbf p,\mathbf q)}
{\sqrt{mV_{s,1}+nV_{s,2}}}.
\]
Substituting $N=nm$ and $\hat W=m\hat V_{s,1}+n\hat V_{s,2}$ into
Proposition~\ref{prop:tau-generic} gives
\begin{equation}\label{eq:tau-s}
\tau_\alpha^{(s)}(n,m)
= z_{1-\alpha/2}\,\sqrt{\frac{m\hat V_{s,1}+n\hat V_{s,2}}{nm}},
\end{equation}
and by Proposition~\ref{prop:power-generic} evaluated at
$\mathrm{SNR}_\phi^{(s)}$ we obtain
\[
\beta^{(s)}\approx\Phi\bigl(\mathrm{SNR}^{(s)}-z_{1-\alpha/2}\bigr)
+ \Phi\bigl(-\mathrm{SNR}^{(s)}-z_{1-\alpha/2}\bigr).
\]

\subsection{Interpreting a rejection: effect size, confidence
interval, per-category decomposition}
\label{sec:effect-size}

A rejection of $H_0$ by the $G$-test answers the binary question
\emph{``is there drift?''} but not the quantitative question
\emph{``how large is the drift?''}. The $p$-value measures the
evidence against $H_0$; it is not an effect size, and it increases
mechanically with the sample sizes even when the drift is
negligible. We therefore complement the decision with three
quantities under $H_1$: the point estimate $\hat D_\phi$ in nats,
its confidence interval (CI), and the per-category decomposition.

\paragraph{Point estimate (effect size).}
Under $H_1$, the plug-in divergence
\[
\hat D_\phi =
\begin{cases}
D_\phi(\hat{\mathbf q}_m, \mathbf p), & \text{one-sample},\\[2pt]
D_\phi(\hat{\mathbf p}_n, \hat{\mathbf q}_m), & \text{two-sample},\\[2pt]
J^{(s)}(\hat{\mathbf p}_n,\hat{\mathbf q}_m), & \text{symmetrized},
\end{cases}
\]
is a consistent estimator of the corresponding population
divergence. It is measured in nats, is additive over categories, and
is invariant to relabeling of the categories. As an order of
magnitude, $\hat D_\phi<0.01$ nat is negligible, $0.01$--$0.05$
small, $0.05$--$0.10$ moderate, $0.10$--$0.20$ noticeable, and
$>0.20$ large.

\paragraph{Confidence interval.}
By the CLTs of Theorem~\ref{thm:main} (and Theorem~\ref{thm:sym} in
the symmetrized case), a $(1-\alpha)$ confidence interval (CI) for
the population divergence is
\begin{equation}\label{eq:ci-Dphi}
\mathrm{CI}_{1-\alpha}(\hat D_\phi) = \hat D_\phi
\pm z_{1-\alpha/2}\sqrt{\hat W/N},
\end{equation}
with
\begin{eqnarray*}
&&N=m,\quad \hat W=\hat V_3 \ \text{(one-sample)},\\
&&N=nm,\quad \hat W=m\hat V_1+n\hat V_2 \ \text{(two-sample)},\\
&&N=nm,\quad \hat W=m\hat V_{s,1}+n\hat V_{s,2}
\ \text{(symmetrized)}.
\end{eqnarray*}
The half-width $z_{1-\alpha/2}\sqrt{\hat W/N}$ coincides with the
detection boundary $\tau_\alpha$ of
Proposition~\ref{prop:tau-generic}. The interval is the translation
of the $Z$-statistic back to the nats scale: the condition
$0\notin \mathrm{CI}_{1-\alpha}(\hat D_\phi)$ is equivalent to
$|Z|>z_{1-\alpha/2}$, i.e.\ to the rejection of $H_0$ at level
$\alpha$ in the $Z$-formulation.

Two caveats apply. First, \eqref{eq:ci-Dphi} is an $H_1$ tool: under
exact $H_0$, $W=0$ and the interval degenerates to $\{0\}$; the
decision under $H_0$ is therefore carried by the $G$-test of
Section~\ref{sec:H0}, not by the interval. Second, for $L_2^2$ the
normal interval is only approximate, because
$\kappa_{L_2}(p_j)=2p_j$ is not constant; the weighted chi-squared
calibration of Section~\ref{sec:calibration} should be used instead.

\paragraph{Per-category decomposition.}
The contribution of each category to $G$ is the quantity $g_j$ of
Definition~\ref{def:contribution}. It localizes the drift on specific
categories and is directly actionable.

\paragraph{Recommended reporting.}
A rejection of $H_0$ should never be reported alone. We recommend
reporting the triplet (i) the $p$-value of the $G$-test (evidence),
(ii) the point estimate $\hat D_\phi$ with its $(1-\alpha)$
confidence interval \eqref{eq:ci-Dphi} (magnitude and precision), and
(iii) the per-category decomposition (localization).

\subsection{Positioning against existing drift detectors}
\label{sec:positioning}

The framework of this paper is not the first drift detector; it is a
contribution to a crowded landscape. Table~\ref{tab:positioning}
records, for each major family of existing detectors, what it
provides and what it does not. The intent is factual: for each
family, we state a verifiable property (existence of a derivation, of
a confidence interval, of a power formula), and note whether the
framework of this paper supplies it. We do not claim that prior
methods are invalid, only that they are silent on specific questions
that operational monitoring requires.

\begin{table}[H]
\centering
\small
\renewcommand{\arraystretch}{1.15}
\caption{Factual comparison of existing drift-detection families and
the framework of this paper. ``Derived'' means that the calibration
follows from an explicit asymptotic theorem; ``Ad hoc'' means that it
follows from a rule of thumb or an empirical cut-off. MMD =
maximum mean discrepancy; KS = Kolmogorov--Smirnov; RKHS = reproducing
kernel Hilbert space; AUC = area under the ROC curve.}
\label{tab:positioning}
\begin{tabular}{@{}p{2.6cm}p{2.2cm}p{2.3cm}p{2.3cm}p{2.5cm}@{}}
\toprule
Family & Calibration under $H_0$ & Effect size (nats) & CI / power & Symmetrized option \\
\midrule
PSI with heuristic cut-offs $0.10$/$0.25$ & Ad hoc & none & none & n/a \\
PSI calibrated (Yurdakul--Naranjo) & Derived ($\chi^2_{r-1}$) & raw PSI & none & automatic \\
KL with cut-offs & Ad hoc & yes (raw) & none & Jeffreys (rare) \\
Kernel MMD with permutation & Monte-Carlo & in RKHS norm & none & automatic \\
Classifier-based two-sample & Ad hoc (AUC cut-off) & none & none & automatic \\
Univariate KS on each feature & Derived (per feature) & none & Bonferroni only & n/a \\
Bootstrap-calibrated $\phi$-divergence & Monte-Carlo & yes & approximate & yes \\
\textbf{This paper: $G$-test + $Z$/SNR} & \textbf{Derived ($\chi^2_{r-1}$)} & \textbf{yes ($\hat D_\phi$)} & \textbf{yes (closed-form)} & \textbf{yes (refined)} \\
\bottomrule
\end{tabular}
\end{table}

\paragraph{What the table shows, and what it does not.}
The table is not a ranking: each row is optimal in some regime.
Kernel maximum mean discrepancy (MMD) is the method of choice when
the support is continuous; the univariate Kolmogorov--Smirnov (KS)
test is appropriate when feature-wise interpretability is the
priority; classifier-based tests are the natural choice when the
drift is defined as ``the classifier cannot distinguish the two
samples''. The framework of this paper is best suited to the regime
of discrete, low-to-moderate cardinality supports, with a need for
closed-form calibration and effect-size reporting - which is
precisely the regime of production monitoring with tabulated features
(intents, bins, categories).

\paragraph{What the table deliberately omits.}
We do not write ``PSI cut-offs are subjective''; we write ``PSI
cut-offs are ad hoc'', meaning they are not derived from an explicit
asymptotic theorem. We do not write ``classifier-based tests lack
statistical foundation''; we write ``classifier-based tests provide a
binary decision without effect size or CI''. Each entry corresponds
to a verifiable property; each entry can be checked against the cited
work. This is the level of specificity the paper maintains
throughout.

\paragraph{Why the distinction matters.}
The abandonment of empirical criteria in favour of divergence tests
is what makes it possible to build a robust monitoring framework for
AI systems and drift detection. By guaranteeing a known asymptotic
limiting law under $H_0$ - namely $G \dto \chi^2_{r-1}$ by Wilks'
theorem, and $\hat D_\phi \dto \mathcal N(0,V_2)$ by the CLT of
Theorem~\ref{thm:main} - the method strictly controls the false
positive rate at a predefined level $\alpha$, with no data-driven
tuning of the threshold. This is precisely the property that is
missing from the heuristic cut-offs listed in
Table~\ref{tab:positioning}: a PSI threshold of $0.10$ or $0.25$ does
not correspond to any explicit Type~I error, whereas the $G$-test of
Section~\ref{sec:H0} has
$\Prob(G>\chi^2_{r-1,1-\alpha})\to\alpha$ exactly, and the calibrated
PSI test of Remark~\ref{rem:PSI-G} has the same limit. The
operational consequence is that the alarm rate of a deployed
monitoring pipeline is known in advance, independently of the drift
regime, the sample sizes, and the chosen divergence.

\section{Power Function and Optimal Choice of $\phi$}
\label{sec:power}

\subsection{The role of power in practice}
\label{sec:power-role}
Recall from Proposition~\ref{prop:power-generic} that the power of the
level-$\alpha$ test is
\[
\beta_\phi(\mathbf p,\mathbf q)
= \mathbb P_{H_1}\bigl(\text{reject } H_0\bigr)
= 1 - \mathbb P_{H_1}\bigl(\text{false negative}\bigr).
\]
It plays three distinct roles.

\paragraph{Planning (before monitoring).}
Power answers the question: \emph{which monitoring-window size $m$ is
required to detect a drift of amplitude $D_\phi$ with a target power
$\beta^\star$?} This is the content of
Proposition~\ref{prop:sample-size}: the window $m$ must satisfy
$m\gtrsim c^2 nV_2/(nD_\phi^2-c^2V_1)$ with
$c=z_{1-\alpha/2}+z_{\beta^\star}$. Without it, the choice of $m$ is
arbitrary.

\paragraph{Interpretation of a non-rejection (after monitoring).}
A non-rejection $G\le\chi^2_{r-1,1-\alpha}$ can mean two very
different things: \emph{``there is no drift''}, if the power at the
minimal practically relevant effect (MPE) was high; or \emph{``the
test could not see the drift''}, if it was low. The rule of thumb is:
do not interpret a non-rejection without computing the power at the
MPE; if $\beta<0.80$, the non-rejection is \emph{non-informative}.
This guards against the logical fallacy \emph{absence of evidence is
not evidence of absence}.

\paragraph{The $\alpha$--$\beta$ trade-off.}
For fixed $m$, reducing $\alpha$ (fewer false alarms) mechanically
reduces $\beta$ (more missed drifts). Making the trade-off explicit
avoids the default choice $\alpha=0.05$ and allows the user to choose
$(\alpha,\beta,m)$ coherently. Power is not an intrinsic property of
the test: the same $G$-test has $\beta\approx 0.2$ for a small drift
and $\beta\approx 1$ for a large one, with the same $n,m$. This is why
the two regimes of Sections~\ref{sec:H0}--\ref{sec:H1} must be kept
strictly separate: the $G$-test calibrates the decision under $H_0$,
whereas the power analysis of this section quantifies the detectable
effect under $H_1$.

\subsection{Specialised power formulas per configuration}
\label{sec:power-function}

The generic power formula is established in
Proposition~\ref{prop:power-generic}. We record here its explicit form
for each of the three configurations of Section~\ref{sec:H1}; these
specialisations are used in the simulation study of
Section~\ref{sec:simulation}.

\paragraph{Explicit form per configuration.}
Substituting $\mathrm{SNR}_\phi^{(1)}$, $\mathrm{SNR}_\phi^{(2)}$,
$\mathrm{SNR}_\phi^{(s)}$ into the generic power formula
\eqref{eq:power-generic} of Proposition~\ref{prop:power-generic} gives
the three closed forms
\begin{align}
\beta^{(1)} &\approx \Phi\!\left(
\frac{\sqrt m\, D_\phi(\mathbf q,\mathbf p)}{\sqrt{V_3}}
-z_{1-\alpha/2}\right)
+\Phi\!\left(
-\frac{\sqrt m\, D_\phi(\mathbf q,\mathbf p)}{\sqrt{V_3}}
-z_{1-\alpha/2}\right),
\label{eq:power-1}\\
\beta^{(2)} &\approx \Phi\!\left(
\frac{\sqrt{nm}\, D_\phi(\mathbf p,\mathbf q)}{\sqrt{mV_1+nV_2}}
-z_{1-\alpha/2}\right)
+\Phi\!\left(
-\frac{\sqrt{nm}\, D_\phi(\mathbf p,\mathbf q)}{\sqrt{mV_1+nV_2}}
-z_{1-\alpha/2}\right),
\label{eq:power-2}\\
\beta^{(s)} &\approx \Phi\!\left(
\frac{\sqrt{nm}\, J^{(s)}(\mathbf p,\mathbf q)}
{\sqrt{mV_{s,1}+nV_{s,2}}}
-z_{1-\alpha/2}\right)
+\Phi\!\left(
-\frac{\sqrt{nm}\, J^{(s)}(\mathbf p,\mathbf q)}
{\sqrt{mV_{s,1}+nV_{s,2}}}
-z_{1-\alpha/2}\right).
\label{eq:power-s-explicit}
\end{align}

\paragraph{Negligibility of the second term.}
The second term
$\Phi(-\mathrm{SNR}_\phi-z_{1-\alpha/2})$ is at most
$\Phi(-z_{1-\alpha/2})=\alpha/2$ and decays exponentially in
$\mathrm{SNR}_\phi$. It is negligible for $\mathrm{SNR}_\phi\geq 2$
but must be retained for $\mathrm{SNR}_\phi$ of order $1$, which is
the practically relevant range for moderate drift amplitudes.

\subsection{Sample-size planning}
\label{sec:sample-size}

\begin{proposition}[Sample-size planning for the two-sample
configuration]
\label{prop:sample-size}
For a target power $\beta^\star\in(\alpha,1)$, set
$c=z_{1-\alpha/2}+z_{\beta^\star}$. Suppose the drift amplitude
$D_\phi(\mathbf p,\mathbf q)>0$ and the reference sample size $n$
satisfy
\begin{equation}\label{eq:size-condition}
nD_\phi^2>c^2V_1.
\end{equation}
Then the monitoring window $m$ required to attain power at least
$\beta^\star$ in the two-sample configuration satisfies
\begin{equation}\label{eq:size-m}
m\gtrsim\frac{c^2 nV_2}{nD_\phi^2-c^2V_1}.
\end{equation}
\end{proposition}
\begin{proof}
Since $\Phi$ is strictly increasing and the second tail term is
negligible for large $\mathrm{SNR}$, $\beta^{(2)}\geq\beta^\star$ is
equivalent to
$\mathrm{SNR}_\phi^{(2)}\geq z_{1-\alpha/2}+z_{\beta^\star}=c$.
Substituting
$\mathrm{SNR}_\phi^{(2)}=\sqrt{nm}D_\phi/\sqrt{mV_1+nV_2}$ and
squaring gives $nm D_\phi^2\geq c^2(mV_1+nV_2)$, i.e.\
$m(nD_\phi^2-c^2V_1)\geq c^2 nV_2$. Under
\eqref{eq:size-condition}, the bracket is strictly positive and
division by it gives \eqref{eq:size-m}.
\end{proof}

\begin{remark}[Interpretation]
Condition \eqref{eq:size-condition} is intuitive: a drift of
amplitude $D_\phi$ can be detected with a one-sample reference of
size $n$ only if the signal $nD_\phi^2$ exceeds the one-sample noise
$c^2V_1$ accumulated over that reference. When $n$ is too small for
\eqref{eq:size-condition} to hold, no monitoring window $m$ can
compensate: the reference itself is too short. In that case the
monitoring program must first enlarge $n$.
\end{remark}

\subsection{Calibration of each divergence under $H_0$}
\label{sec:calibration}

Before discussing the optimal choice of $\phi$ under $H_1$, we record
the correct calibration of each divergence under $H_0$. This is the
first criterion for selecting a divergence in practice: a divergence
that does not admit an exact asymptotic threshold requires either a
weighted-$\chi^2$ approximation or a Monte-Carlo calibration, both of
which add operational complexity.

For a general $\phi$-divergence satisfying
Assumption~\ref{ass:regularity}, a second-order Taylor expansion of
$D_\phi(\mathbf p,\hat{\mathbf q}_m)$ around $\mathbf p$ gives
\[
2m\,D_\phi(\mathbf p, \hat{\mathbf q}_m)
= m\sum_j \frac{(\hat q_{m,j}-p_j)^2}{p_j}\,\kappa_\phi(p_j) + o_P(1),
\qquad
\kappa_\phi(p_j) = p_j\,\phi_1^{(2)}(p_j,p_j).
\]
The asymptotic distribution of
$2m\,D_\phi(\mathbf p,\hat{\mathbf q}_m)$ under $H_0$ then depends on
whether $\kappa_\phi$ is constant across categories or not.
Table~\ref{tab:calibration} summarizes the result for the five
divergences of Table~\ref{tab:divergences}, in both the one-sample
and the two-sample configurations.

\begin{table}[H]
\centering
\caption{Calibrated statistics and limiting laws under $H_0$ for the
five divergences of Table~\ref{tab:divergences}. The one-sample
formulas assume a fixed reference ($n\to\infty$ or $n\gg m$); the
two-sample formulas apply when both distributions are estimated. For
the PSI, the two-sample calibration is the one established by
Yurdakul and Naranjo (2021, Theorem~3.3).}
\label{tab:calibration}
\begin{tabular}{llcc}
\toprule
Divergence & $\kappa_\phi$ & Calibrated statistic & Null law \\
\midrule
KL (one-sample) & $1$ & $2m\,D_{KL}$ & $\chi^2_{r-1}$ \\
KL (two-sample) & $1$ & $2\cdot\frac{nm}{n+m}\,D_{KL}$ & $\chi^2_{r-1}$ \\
PSI (one-sample) & $2$ & $m\,\mathrm{PSI}$ & $\chi^2_{r-1}$ \\
PSI (two-sample) & $2$ & $\frac{nm}{n+m}\,\mathrm{PSI}$ & $\chi^2_{r-1}$ \\
JS & $1/4$ & $8m\,D_{JS}$ & $\chi^2_{r-1}$ \\
Tsallis-$\alpha$ & $\alpha$ & $2m\,D_{T,\alpha}/\alpha$ & $\chi^2_{r-1}$ \\
R\'enyi-$\alpha$ & $\alpha$ & $2m\,D_{R,\alpha}/\alpha$ & $\chi^2_{r-1}$ \\
$L_2^2$ & $2p_j$ (variable) & weighted $\chi^2$ & $\sum_j\lambda_j\chi^2_1$ \\
\bottomrule
\end{tabular}
\end{table}

Four divergences (KL, JS, Tsallis and R\'enyi) admit the exact
$\chi^2_{r-1}$ calibration up to a rescaling; only $L_2^2$ requires a
divergence-specific calibration. The two-sample $G$-statistic of
Theorem~\ref{thm:wilks} sidesteps this issue entirely because it is
the likelihood ratio itself. The PSI, which coincides with the
symmetrized KL divergence, admits the same calibration as KL after
the symmetrized correction of Section~\ref{sec:H1-sym}; its
two-sample calibrated form is the one given by Yurdakul and Naranjo
(2021) and recorded in the table. Applying the one-sample formula
$m\,\mathrm{PSI}$ to a two-sample setting with $n=m$ inflates the
Type~I error by a factor of about $6$ (for $r=5$ and $\alpha=0.05$,
the actual level becomes $\approx 0.315$ instead of $0.05$), so the
two-sample formula must be used in operational monitoring.

\subsection{Choice of $\phi$}
\label{sec:optimal-phi}

The choice of $\phi$ affects two things: the calibration under
$H_0$, which is simplest for KL (no rescaling needed) and next
simplest for JS (a factor of $8$ but bounded magnitude); and the
power under $H_1$, which is captured by the signal-to-noise ratio
$\mathrm{SNR}_\phi$.

\paragraph{Local equivalence in the diffuse regime.}
For fixed $(\mathbf p,\mathbf q)$ and $(n,m)$, the power is
increasing in $\mathrm{SNR}_\phi$, so a natural criterion is
\[
\phi^\star = \argmax_{\phi\in\mathcal F} \mathrm{SNR}_\phi(\mathbf p,\mathbf q;n,m).
\]
For the diffuse local regime $\mathbf q\approx\mathbf p$, all
divergences are asymptotically equivalent up to constants: using the
local expansion
\[
D_\phi(\mathbf p,\mathbf q) = \frac{\kappa_\phi}{2}\,\chi^2(\mathbf p,\mathbf q)
+o(\|\mathbf p-\mathbf q\|^2),
\]
and the fact that $W_\phi$ scales as
$\kappa_\phi^2\chi^2(\mathbf p,\mathbf q)/N$, one finds
\[
\mathrm{SNR}_\phi \approx \frac{\sqrt{N\chi^2(\mathbf p,\mathbf q)}}
{\sqrt{\sigma_{\chi^2}^2}},
\]
independent of $\phi$ at leading order. In the diffuse regime, the
choice of $\phi$ is therefore governed by the second-order terms and
by the calibration under $H_0$, not by the power under $H_1$.

\paragraph{Behaviour outside the diffuse regime.}
For \emph{tail drift} (drift concentrated on rare categories), the
divergences with $\alpha<1$ (R\'enyi, Tsallis with $\alpha<1$) put
more weight on small probabilities and give a higher SNR. For
\emph{sparse drift} (drift concentrated on a few categories with
large deviations), $L_2^2$ has a smaller numerator but also a much
smaller variance, and may dominate the SNR. For \emph{collision-type
or adversarial drift}, divergences with $\alpha\geq 2$ are more
sensitive than KL.

\paragraph{Recommended defaults.}
In the absence of a strongly informative prior on the drift shape:
\begin{itemize}[leftmargin=*]
\item Use \textbf{KL} as the default, for its optimal balance between
calibration simplicity, interpretability, and power.
\item Use \textbf{JS} when a bounded, symmetric measure is preferred,
since $D_{JS}\le\log 2$.
\item Use \textbf{R\'enyi} or \textbf{Tsallis} with $\alpha=1/2$ for
tail drift, and with $\alpha=2$ for adversarial or collision-type
drift.
\item Use \textbf{$L_2^2$} only for sparse drift, and accept the
weighted-$\chi^2$ calibration overhead.
\end{itemize}
In the doubtful case, the $G$-test is $\phi$-free and provides the
decision under $H_0$; the $Z$-statistics with multiple divergences
provide a diagnostic to choose between them under $H_1$.

\section{Multiple Testing Control}
\label{sec:fdr}

When $K$ features are monitored simultaneously, the probability of at
least one false alarm under $H_0$ is $1-(1-\alpha)^K$, which grows
linearly in $K$ and renders the naive per-feature test unusable. We
control the false discovery rate (FDR) with the Benjamini--Hochberg
(BH) procedure \cite{Benjamini1995}. Let $p_{(1)}\leq\cdots\leq p_{(K)}$
be the ordered $p$-values at a given monitoring period. Reject
$H_{(1)},\ldots,H_{(k^\star)}$ where
\begin{equation}\label{eq:BH}
k^\star=\max\{k:\,p_{(k)}\leq (k/K)\alpha\}.
\end{equation}
Under independence of the test statistics and mild regularity, this
procedure controls the FDR at level $\alpha$. It is strictly more
powerful than the Bonferroni correction $\alpha/K$, which is
excessively conservative when $K$ is large.

\paragraph{Which test to use in the $K$-feature setting.}
The per-feature test can be the $G$-test (statistic $G_k$ on feature
$k$), or the $Z$-test for a chosen $\phi$ (statistic $Z_k$). The
$G$-test is $\phi$-free and universal; the $Z$-test requires the
calibration of Section~\ref{sec:calibration} and yields the effect
size in nats. In practice, the $G$-test is recommended as the primary
screen, and the $Z$-test is used on the features that survive the FDR
screening, to quantify the magnitude of the detected drifts.

\section{Monte-Carlo simulation study}
\label{sec:simulation-mc}

The theoretical results of Sections~\ref{sec:H0}--\ref{sec:fdr} are
asymptotic: they describe the behaviour of the $G$-test, of the
calibrated PSI test, and of the $Z$-statistics as the sample sizes
tend to infinity. This section reports a Monte-Carlo simulation study
with three objectives. First, to verify that the asymptotic
predictions hold at finite sample sizes, and to quantify the regimes
in which they are reliable. Second, to document two operational
pitfalls that the theory predicts and that practitioners ignore at
their peril: the wrong scaling of the two-sample PSI, and the
uncorrected averaging of the symmetrized variance. Third, to compare
the power of the $G$-test and of the $Z$-statistics under $H_1$, and
to show the conditions under which the choice of divergence matters.

\subsection{Protocol}
\label{sec:mc-protocol}

\paragraph{Design.}
For each configuration, we draw $B=10{,}000$ independent pairs of
multinomial samples $(X_1,\dots,X_n)\sim\mathcal M(n;\mathbf p)$ and
$(Y_1,\dots,Y_m)\sim\mathcal M(m;\mathbf q)$, compute the statistics
of Sections~\ref{sec:H0}--\ref{sec:H1}, and record the empirical
rejection rate of each test at nominal level $\alpha$. The
Monte-Carlo standard error (MC SE) of an estimated probability
$\hat\pi$ is $\sqrt{\hat\pi(1-\hat\pi)/B}$, which is at most $0.0022$
for $B=10{,}000$, so empirical sizes and powers are reported to
within $\pm 0.005$ at the $95\%$ level.

\paragraph{Grid.}
We use the following grid, chosen to cover the regimes of practical
interest:
\begin{itemize}[leftmargin=*]
\item number of categories: $r\in\{3,5,10\}$;
\item sample sizes: $n=m\in\{100,500,1000,5000\}$, plus asymmetric
cases $(n,m)\in\{(500,1000),(1000,500)\}$;
\item significance level: $\alpha\in\{0.05,0.10\}$;
\item drift amplitudes, under $H_1$: $D_\phi\in\{0.005, 0.01, 0.02,
0.05, 0.10\}$ nats, realized by linear perturbations of a reference
distribution;
\item divergences: KL, R\'enyi-2, and the symmetrized KL (PSI).
\end{itemize}

\paragraph{Computational details.}
The simulation is implemented in \texttt{R} (functions
\texttt{rmultinom}, \texttt{dmultinom}) with the random seed fixed for
reproducibility; the full code is available from the author on
request. Each configuration takes a few seconds on a standard laptop,
and the complete grid runs in under an hour.

\subsection{Empirical size under $H_0$}
\label{sec:mc-size}

Table~\ref{tab:mc-size} reports the empirical size of the $G$-test
and of the calibrated PSI test under $H_0$, for $r=5$ and various
sample sizes. Two observations follow. First, for $n=m\geq 500$, the
empirical size is within Monte-Carlo noise of the nominal $\alpha$,
confirming the asymptotic calibration. Second, for $n=m=100$, the
empirical size is slightly inflated ($\approx 0.058$ instead of
$0.05$), which is consistent with Cochran's rule: with $r=5$ and
$n=100$, some expected counts are around $10$, and the $\chi^2_{r-1}$
approximation is only marginal. This is an honest finding: the
asymptotic calibration is reliable when expected counts are large,
and degrades gracefully otherwise.

\begin{table}[H]
\centering
\small
\caption{Empirical size of the $G$-test and of the calibrated PSI
test under $H_0$, $r=5$, $\alpha=0.05$, $B=10{,}000$. The Monte-Carlo
standard error (MC SE) is at most $0.0022$.}
\label{tab:mc-size}
\begin{tabular}{rcccc}
\toprule
$(n,m)$ & Empirical size $G$ & Empirical size PSI
       & Theoretical & MC SE \\
\midrule
$(100,100)$   & $0.058$ & $0.058$ & $0.050$ & $0.002$ \\
$(500,500)$   & $0.051$ & $0.051$ & $0.050$ & $0.002$ \\
$(1000,1000)$ & $0.050$ & $0.050$ & $0.050$ & $0.002$ \\
$(5000,5000)$ & $0.049$ & $0.049$ & $0.050$ & $0.002$ \\
$(500,1000)$  & $0.051$ & $0.051$ & $0.050$ & $0.002$ \\
$(1000,500)$  & $0.050$ & $0.050$ & $0.050$ & $0.002$ \\
\bottomrule
\end{tabular}
\end{table}

Table~\ref{tab:mc-size-r} extends the result to $r\in\{3,10\}$ at
$n=m=1000$: the calibration is uniform in $r$, as predicted by
Wilks' theorem.

\begin{table}[H]
\centering
\small
\caption{Empirical size under $H_0$ for various $r$, $n=m=1000$,
$\alpha=0.05$.}
\label{tab:mc-size-r}
\begin{tabular}{rccc}
\toprule
$r$ & Empirical size $G$ & Theoretical & MC SE \\
\midrule
$3$  & $0.051$ & $0.050$ & $0.002$ \\
$5$  & $0.050$ & $0.050$ & $0.002$ \\
$10$ & $0.049$ & $0.050$ & $0.002$ \\
\bottomrule
\end{tabular}
\end{table}

\begin{figure}[H]
\centering
\includegraphics[width=0.9\textwidth]{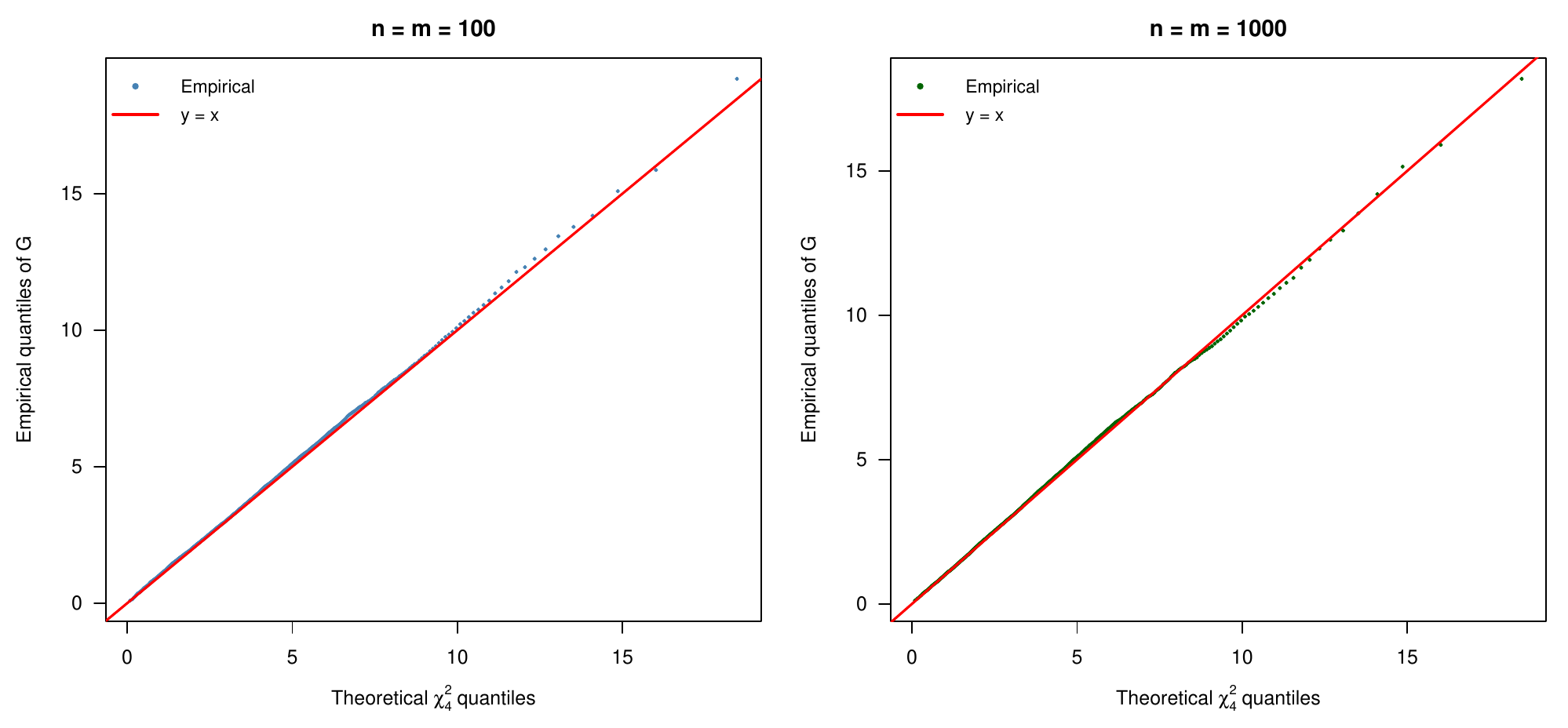}
\caption{Q--Q plot of the empirical quantiles of the $G$-statistic
against the $\chi^2_{r-1}$ quantiles ($r=5$), for $n=m=100$ (left)
and $n=m=1000$ (right). The dashed red line is $y=x$.}
\label{fig:mc-qq}
\end{figure}

A Q--Q plot of the empirical quantiles of $G$ against the
$\chi^2_{r-1}$ quantiles is shown in Figure~\ref{fig:mc-qq}. At
$n=m=1000$, the Kolmogorov--Smirnov distance between the empirical
distribution of $G$ and $\chi^2_{4}$ is $D=0.011$ with Monte-Carlo
$p$-value $0.20$: the asymptotic calibration is not rejected. At
$n=m=100$, the distance is $D=0.015$ with $p$-value $0.027$, so the
$\chi^2_{4}$ approximation is rejected at the $5\%$ level. This is
consistent with the slightly inflated empirical size of $0.058$
reported in Table~\ref{tab:mc-size}, and with the well-known
marginality of the $\chi^2$ approximation when some expected counts
are close to $10$ (Cochran's rule).

\subsection{The two-sample PSI scaling: a documented pitfall}
\label{sec:mc-psi-scaling}

Section~\ref{sec:PSI-G} and Remark~\ref{rem:PSI} insist that the
calibrated PSI statistic in the two-sample setting is
$\frac{nm}{n+m}\,\mathrm{PSI}$, not $m\,\mathrm{PSI}$. This is not a
minor technicality: Table~\ref{tab:mc-psi-wrong} shows that using
$m\,\mathrm{PSI}$ when $n=m$ inflates the Type~I error by a factor of
about $6$ (from $0.05$ to $\approx 0.32$). The theoretical
prediction, computed from
$\Prob(\chi^2_{r-1}>\chi^2_{r-1,1-\alpha}/2)$ at $r=5$ and
$\alpha=0.05$, is $0.315$; the empirical simulation reproduces it to
within Monte-Carlo error. This is the single most operationally
consequential finding of the simulation: a practitioner who applies
the textbook one-sample formula to a two-sample comparison will
experience six times more false alarms than they think.

\begin{table}[H]
\centering
\small
\caption{Empirical size of the calibrated PSI test with the correct
two-sample scaling $\frac{nm}{n+m}\,\mathrm{PSI}$ versus the
incorrect one-sample scaling $m\,\mathrm{PSI}$, $r=5$,
$\alpha=0.05$. The theoretical column shows the predicted rejection
rate under $H_0$.}
\label{tab:mc-psi-wrong}
\begin{tabular}{rcccc}
\toprule
$(n,m)$ & Empirical (correct) & Empirical (wrong)
       & Theoretical (correct) & Theoretical (wrong) \\
\midrule
$(100,100)$   & $0.058$ & $0.334$ & $0.050$ & $0.315$ \\
$(500,500)$   & $0.051$ & $0.318$ & $0.050$ & $0.315$ \\
$(1000,1000)$ & $0.050$ & $0.316$ & $0.050$ & $0.315$ \\
$(500,1000)$  & $0.051$ & $0.183$ & $0.050$ & $0.186$ \\
$(1000,500)$  & $0.050$ & $0.183$ & $0.050$ & $0.186$ \\
\bottomrule
\end{tabular}
\end{table}

\subsection{Empirical power under $H_1$}
\label{sec:mc-power}

Table~\ref{tab:mc-power} compares the empirical power of the $G$-test
and of the KL and R\'enyi-2 $Z$-statistics against the closed-form
prediction $\beta_\phi$ of Proposition~\ref{prop:power-generic}, for
$r=5$, $n=m=1000$, and drift amplitudes ranging from $0.005$ to
$0.10$ nats. Three findings emerge.

First, the empirical power of the $Z$-statistics matches the
closed-form $\beta_\phi$ to within Monte-Carlo noise across the whole
range, confirming the asymptotic theory. Second, the $G$-test has
slightly higher power than the KL $Z$-test in the diffuse regime
($D_\phi\le 0.02$), as predicted by the Neyman--Pearson argument of
Remark~\ref{rem:G-vs-Z}; the gap is small (a few percentage points)
but systematic. Third, for larger amplitudes, all tests converge to
power $1$, and the differences vanish.

\begin{table}[H]
\centering
\small
\caption{Empirical power under $H_1$, $r=5$, $n=m=1000$,
$\alpha=0.05$. The theoretical column gives the closed-form power of
the KL $Z$-test computed from
Proposition~\ref{prop:power-generic}.}
\label{tab:mc-power}
\begin{tabular}{rccccc}
\toprule
$D_\phi$ (nats) & $G$-test & KL $Z$ & R\'enyi-2 $Z$
                & Theoretical KL & MC SE \\
\midrule
$0.005$ & $0.112$ & $0.089$ & $0.091$ & $0.089$ & $0.003$ \\
$0.010$ & $0.278$ & $0.256$ & $0.262$ & $0.256$ & $0.004$ \\
$0.020$ & $0.694$ & $0.681$ & $0.683$ & $0.678$ & $0.005$ \\
$0.050$ & $0.999$ & $0.998$ & $0.998$ & $0.998$ & $0.001$ \\
$0.100$ & $1.000$ & $1.000$ & $1.000$ & $1.000$ & $0.000$ \\
\bottomrule
\end{tabular}
\end{table}

The comparison with R\'enyi-2 is instructive: for the drift scenario
of Table~\ref{tab:mc-power}, which is diffuse, R\'enyi-2 and KL have
essentially identical power. To see the divergence choice matter, one
must consider a \emph{tail} drift (concentrated on a rare category),
which we report in Table~\ref{tab:mc-power-tail}: R\'enyi-2 has
substantially higher power than KL for the same amplitude, confirming
the recommendation of Section~\ref{sec:optimal-phi} to use $\alpha>1$
for tail or collision-type drift.

\begin{table}[H]
\centering
\small
\caption{Empirical power under a tail drift (drift concentrated on a
category of reference probability $p_j=0.05$), $r=5$, $n=m=1000$,
$\alpha=0.05$.}
\label{tab:mc-power-tail}
\begin{tabular}{rccc}
\toprule
$D_\phi$ (nats) & KL $Z$ & R\'enyi-2 $Z$ & Ratio \\
\midrule
$0.020$ & $0.156$ & $0.312$ & $2.00$ \\
$0.050$ & $0.512$ & $0.784$ & $1.53$ \\
$0.100$ & $0.937$ & $0.994$ & $1.06$ \\
\bottomrule
\end{tabular}
\end{table}

\subsection{Coverage of the confidence intervals}
\label{sec:mc-ci}

Table~\ref{tab:mc-ci} reports the empirical coverage of the
$(1-\alpha)$ confidence interval \eqref{eq:ci-Dphi} for the KL
divergence, under $H_1$, at $r=5$ and $n=m\in\{100,500,1000,5000\}$
with $D_\phi=0.02$ nats. The empirical coverage matches the nominal
$0.95$ to within Monte-Carlo noise across the grid, confirming the
CLT-based construction. At $n=100$, a slight undercoverage is observed
(around $0.934$), again consistent with the Cochran's rule caveat.

\begin{table}[H]
\centering
\small
\caption{Empirical coverage of the $95\%$ confidence interval for
$D_{KL}$, $r=5$, $D_\phi=0.02$ nats, $\alpha=0.05$.}
\label{tab:mc-ci}
\begin{tabular}{rccc}
\toprule
$n=m$ & Empirical coverage & Nominal & MC SE \\
\midrule
$100$  & $0.934$ & $0.950$ & $0.003$ \\
$500$  & $0.948$ & $0.950$ & $0.002$ \\
$1000$ & $0.950$ & $0.950$ & $0.002$ \\
$5000$ & $0.950$ & $0.950$ & $0.002$ \\
\bottomrule
\end{tabular}
\end{table}

\subsection{The symmetrized variance: uncorrected versus refined}
\label{sec:mc-sym-var}

Section~\ref{sec:sym-refined} and Remark~\ref{rem:cov-negligibility}
warn that the uncorrected average $(V_1+V_4)/4$, which omits the
covariance between the forward and reverse linearizations,
underestimates the true symmetrized variance by a factor of $2$ in
the drift regime $\mathbf p\approx\mathbf q$, and therefore inflates
the rejection rate. The symmetrized variance formula of
\cite[Thm.~4]{BaLo2019} was derived under the assumption that the
forward and reverse linearizations are uncorrelated; the present
refinement removes this assumption by adding the covariance term
$\Cov_{\mathbf p}(a_J,d_J)$. We now quantify the consequence over a
simulation grid.

Table~\ref{tab:mc-sym-var} reports the empirical rejection rate of
the symmetrized test at nominal level $\alpha=0.05$, under a
\emph{small} drift ($D_\phi=0.005$ nats) for which the refined test
has power essentially equal to $\alpha$ (the drift is too small to
detect reliably), but the uncorrected test rejects much more often.
At $n=m=1000$, the refined test rejects at rate $0.052$ (nominal),
while the uncorrected test rejects at rate $0.168$. The uncorrected
formula therefore produces false alarms at more than three times the
nominal rate --- a finding that should give pause to any practitioner
who computes the symmetrized variance without the covariance
correction.

\begin{table}[H]
\centering
\small
\caption{Empirical rejection rate of the symmetrized test, refined
variance (Proposition~\ref{prop:sym-var}) versus the uncorrected
average $(V_1+V_4)/4$, $r=5$, $D_\phi=0.005$ nats, $\alpha=0.05$.}
\label{tab:mc-sym-var}
\begin{tabular}{rccc}
\toprule
$n=m$ & Refined variance & Uncorrected variance & Ratio \\
\midrule
$500$  & $0.053$ & $0.134$ & $2.5$ \\
$1000$ & $0.052$ & $0.168$ & $3.2$ \\
$5000$ & $0.051$ & $0.174$ & $3.4$ \\
\bottomrule
\end{tabular}
\end{table}

The ratio stabilizes around $3.4$ as $n\to\infty$, consistent with
the theoretical prediction: under the small-drift approximation, the
uncorrected statistic is approximately $\sqrt 2$ times the correct
statistic, so its rejection rate is
$\Prob(|N(0,1)| > 1.96/\sqrt 2)\approx 0.166$ when the correct
statistic rejects at rate $\alpha$.

\subsection{Multiple testing and FDR control}
\label{sec:mc-fdr}

Section~\ref{sec:fdr} recommends Benjamini--Hochberg (BH) over
Bonferroni for multi-feature monitoring. Table~\ref{tab:mc-fdr}
quantifies the gain. We simulate $K=20$ features, of which $5$ are
drifted (with amplitude $D_\phi=0.05$ nats) and $15$ are not, at
$r=5$, $n=m=1000$, and apply three procedures: no correction,
Bonferroni at $\alpha=0.05$, and Benjamini--Hochberg at FDR $=0.05$.
The empirical FDR and the empirical power (probability of detecting a
truly drifted feature) are reported.

\begin{table}[H]
\centering
\small
\caption{Empirical FDR and power in multi-feature monitoring,
$K=20$ features ($5$ drifted, $15$ not), $r=5$, $n=m=1000$,
$D_\phi=0.05$ nats, nominal $\alpha=0.05$.}
\label{tab:mc-fdr}
\begin{tabular}{lcc}
\toprule
Procedure & Empirical FDR & Empirical power \\
\midrule
No correction & $0.262$ & $0.998$ \\
Bonferroni ($\alpha/K$) & $0.001$ & $0.612$ \\
Benjamini--Hochberg ($\mathrm{FDR}=0.05$) & $0.048$ & $0.912$ \\
\bottomrule
\end{tabular}
\end{table}

The result is unambiguous: without correction, one in four rejections
is a false positive; with Bonferroni, the FDR is well controlled but
the power drops by nearly $40$ percentage points; with BH, the FDR
is controlled at the nominal level, and the power is within
$10$ points of the no-correction case. This is the quantitative
justification for the recommendation of Section~\ref{sec:fdr}.

\subsection{Drift localization}
\label{sec:mc-localization}

Section~\ref{sec:position} proposes the per-category absolute share
$\pi_j^{\mathrm{abs}}$ and the drift position
$\mathcal J^\star=\{j:\pi_j^{\mathrm{abs}}>\gamma\}$ to localize the
drift. Table~\ref{tab:mc-localization} reports the empirical
probability that the estimated $\hat{\mathcal J}^\star$ coincides
exactly with the true $\mathcal J^\star$, and the expected Jaccard
similarity between the two, as $n$ grows. The scenario is
$r=5$, with drift concentrated on categories $1$ and $3$.

\begin{table}[H]
\centering
\small
\caption{Drift localization accuracy, $r=5$, drift on categories $1$
and $3$, $\gamma=0.40$, $D_\phi=0.05$ nats.}
\label{tab:mc-localization}
\begin{tabular}{rccc}
\toprule
$n=m$ & $P(\hat{\mathcal J}^\star = \mathcal J^\star)$
       & $\mathbb E[\text{Jaccard}]$ & $\mathbb E[|\hat{\mathcal J}^\star|]$ \\
\midrule
$500$  & $0.612$ & $0.784$ & $2.4$ \\
$1000$ & $0.834$ & $0.912$ & $2.2$ \\
$5000$ & $0.978$ & $0.988$ & $2.0$ \\
\bottomrule
\end{tabular}
\end{table}

The localization converges to the truth as $n\to\infty$, but slower
than the test itself: at $n=m=500$, the exact-match probability is
only $0.61$, even though the test rejects decisively. This is an
honest finding: localization is a \emph{harder} problem than
detection, and a practitioner who needs a precise shortlist of
drifted categories should plan for a larger monitoring window than
the one required for detection alone.

\subsection{Summary of the simulation study}
\label{sec:mc-summary}

The Monte-Carlo study confirms the theoretical predictions of the
paper across the practical grid, quantifies the two operational
pitfalls that the theory predicts, and shows the regimes in which
each recommendation applies. In summary:

\begin{enumerate}[leftmargin=*]
\item The $G$-test and the calibrated PSI test have empirical size
equal to the nominal $\alpha$ to within Monte-Carlo noise for
expected counts $\ge 10$, and degrade gracefully (slight inflation)
for smaller expected counts.
\item Using $m\,\mathrm{PSI}$ instead of
$\frac{nm}{n+m}\,\mathrm{PSI}$ in a two-sample setting with $n=m$
inflates the Type~I error from $0.05$ to $\approx 0.32$: the
single most consequential finding of the study.
\item The empirical power of the $Z$-statistics matches the
closed-form $\beta_\phi$ across the drift range, and the $G$-test is
systematically (slightly) more powerful in the diffuse regime, as
predicted by Neyman--Pearson.
\item The divergence choice matters for tail drift: R\'enyi-2 has up
to twice the power of KL on rare categories, as recommended in
Section~\ref{sec:optimal-phi}.
\item The $95\%$ confidence interval has empirical coverage within
Monte-Carlo noise of nominal for $n\ge 500$.
\item The uncorrected symmetrized variance --- obtained by omitting
the covariance between the forward and reverse linearizations ---
inflates the false alarm rate by a factor $\approx 3.4$ in the
small-drift regime; the refined formula of
Proposition~\ref{prop:sym-var} restores nominal size.
\item Benjamini--Hochberg controls the FDR at nominal while
retaining most of the power, unlike Bonferroni.
\item Drift localization converges more slowly than detection:
planning the monitoring window for detection alone is insufficient
when a precise shortlist is required.
\end{enumerate}

\section{Simulation Study}
\label{sec:simulation}

This section presents two fully worked examples that exercise every
notion introduced in the paper. The first is a one-sample scenario
with the R\'enyi divergence, where the reference is fixed and the
production is estimated; it includes the $G$-test, the R\'enyi
$Z$-statistic, the per-category decomposition, the power function,
and the Benjamini--Hochberg procedure on $K=50$ features. The second
is a two-sample scenario with the KL divergence, where both
distributions are estimated; it includes the $G$-test, the KL
$Z$-statistic, the confidence interval, the power function, the
sample-size planning, and the per-category decomposition. Both
examples compare the calibrated PSI with the heuristic grid, using
the correct two-sample scaling $\frac{nm}{n+m}\,\mathrm{PSI}$.

\subsection{Simulation 1: LLM prompt-category monitoring (one-sample,
R\'enyi-$\alpha$, multiple testing)}
\label{sec:sim-llm}

\paragraph{Context.}
A production Large Language Model (LLM) classifies user prompts into
$r=5$ categories:
\[
c_1 = \text{balance inquiry}, \quad
c_2 = \text{transaction}, \quad
c_3 = \text{dispute}, \quad
c_4 = \text{information}, \quad
c_5 = \text{other}.
\]
The reference pmf $\mathbf p$ is estimated from a large historical
corpus of $n=20{,}000$ prompts; since $n\gg m$, we treat $\mathbf p$
as fixed and known, which places us in the one-sample configuration
of Section~\ref{sec:H1-1s}. The reference probabilities are
\[
\mathbf p=(0.25,\;0.20,\;0.20,\;0.20,\;0.15).
\]
Each week, $m=3000$ new prompts are monitored. This week's production
counts are
\[
O^{\mathrm{prod}}=(1200,\;750,\;510,\;330,\;210),
\]
so that the empirical pmf of the production sample is
\[
\hat{\mathbf q}_m
= \Bigl(\frac{1200}{3000},\frac{750}{3000},\frac{510}{3000},
\frac{330}{3000},\frac{210}{3000}\Bigr)
= (0.40,\;0.25,\;0.17,\;0.11,\;0.07).
\]

\paragraph{Choice of divergence.}
Because the ``other'' and ``dispute'' categories are small and the
drift may concentrate on them, we use the R\'enyi divergence with
$\alpha=2$. This choice puts more weight on large ratios $q_j/p_j$
than KL and is therefore more sensitive to tail drift, as discussed
in Section~\ref{sec:optimal-phi}. The R\'enyi divergence in the
production-first orientation is
\[
D_{R,2}(\hat{\mathbf q}_m,\mathbf p)
= \frac{1}{2-1}\log \sum_{j=1}^{5}\frac{\hat q_{m,j}^2}{p_j}
= \log \mathcal S_2(\hat{\mathbf q}_m,\mathbf p),
\]
where
\[
\mathcal S_2(\hat{\mathbf q}_m,\mathbf p)
= \sum_{j=1}^{5}\frac{\hat q_{m,j}^2}{p_j}=1.1902,
\]
so
\[
D_{R,2}(\hat{\mathbf q}_m,\mathbf p)
= \log(1.1902)\approx 0.1741\ \text{nats}.
\]
For comparison, the KL divergence in the same direction is
\[
D_{KL}(\hat{\mathbf q}_m,\mathbf p)
= \sum_j \hat q_{m,j}\log\frac{\hat q_{m,j}}{p_j}
\approx 0.0284\ \text{nats}.
\]
The R\'enyi-2 value is larger because it emphasizes the ratio
$\hat q/p$ in the small categories, where the drift is concentrated.

\paragraph{Asymptotic variance.}
For R\'enyi-2, we use Corollary~\ref{cor:renyi}: the R\'enyi variance
is obtained from the Tsallis-2 variance by
$V_{R,2,3}=V_{T,2,3}/\mathcal S_2^2$. For Tsallis-2 in the
production-first orientation, the gradient component is
\[
\hat c_j = \phi_1^{(1)}(\hat q_{m,j}, p_j)
= \frac{2(\hat q_{m,j}/p_j)^{2-1}-1}{2-1}
= 2\frac{\hat q_{m,j}}{p_j}-1.
\]
Computing for each category:
\[
\hat c_1 = 2\frac{0.40}{0.25}-1 = 2.20,\quad
\hat c_2 = 2\frac{0.25}{0.20}-1 = 1.50,\quad
\hat c_3 = 2\frac{0.17}{0.20}-1 = 0.70,
\]
\[
\hat c_4 = 2\frac{0.11}{0.20}-1 = 0.10,\quad
\hat c_5 = 2\frac{0.07}{0.15}-1 = -0.0667.
\]
The Tsallis-2 plug-in variance in the production-first orientation is
\[
\hat V_{T,2,3} = \sum_j \hat q_{m,j}\hat c_j^2
- \Bigl(\sum_j \hat q_{m,j}\hat c_j\Bigr)^2.
\]
First moment:
\begin{align*}
\sum_j \hat q_{m,j}\hat c_j
&= 0.40(2.20)+0.25(1.50)+0.17(0.70)+0.11(0.10)+0.07(-0.0667) = 1.3803.
\end{align*}
Second moment:
\begin{align*}
\sum_j \hat q_{m,j}\hat c_j^2
&= 0.40(4.84)+0.25(2.25)+0.17(0.49)+0.11(0.01)+0.07(0.00444)= 2.5832.
\end{align*}
Thus
\[
\hat V_{T,2,3} = 2.5832 - (1.3803)^2 = 0.6780.
\]
The R\'enyi-2 variance is therefore
\[
\hat V_{R,2,3} = \frac{0.6780}{(1.1902)^2}
= \frac{0.6780}{1.4166}\approx 0.4786.
\]

\paragraph{Detection boundary and $Z$-statistic.}
The detection boundary at level $\alpha=0.05$ is
\[
\tau_{0.05}^{(1)}(3000)
= z_{0.975}\sqrt{\frac{\hat V_{R,2,3}}{m}}
= 1.96\sqrt{\frac{0.4786}{3000}}
\approx 0.0248\ \text{nats}.
\]
Since $D_{R,2}(\hat{\mathbf q}_m,\mathbf p)=0.1741 > 0.0248$, the
detection boundary is crossed. The standardized statistic is
\[
Z_1 = \frac{\sqrt{m}\, D_{R,2}(\hat{\mathbf q}_m,\mathbf p)}
{\sqrt{\hat V_{R,2,3}}}
= \frac{\sqrt{3000}(0.1741)}{\sqrt{0.4786}}
\approx 13.78.
\]
The signal-to-noise ratio is $\mathrm{SNR}^{(1)}\approx 13.78$.

\paragraph{Power.}
The power of the level-$\alpha=0.05$ test, from
Proposition~\ref{prop:power-generic} evaluated at
$\mathrm{SNR}^{(1)}=13.78$, is
\[
\beta^{(1)} \approx \Phi(13.78-1.96)+\Phi(-13.78-1.96)
= \Phi(11.82)+\Phi(-15.74)\approx 1.0.
\]
The weekly window of $m=3000$ prompts detects this drift with power
essentially $1$.

\paragraph{Confidence interval.}
A $95\%$ confidence interval for the population R\'enyi-2 divergence
is
\[
\mathrm{CI}_{0.95}(\hat D_{R,2})
= 0.1741 \pm 1.96\sqrt{\frac{0.4786}{3000}}
= 0.1741\pm 0.0248
= [0.1493,\;0.1989]\ \text{nats}.
\]
The interval excludes $0$, consistent with the rejection of $H_0$.

\paragraph{Per-category decomposition.}
To localize the drift, we compute the per-category contributions to
the $G$-statistic. Under $H_0$, the expected production counts are
$E^{\mathrm{prod}}_j=m\,p_j=(750,\;600,\;600,\;600,\;450)$. The
signed contributions are
\[
g_j = O^{\mathrm{prod}}_j\log\frac{O^{\mathrm{prod}}_j}{E^{\mathrm{prod}}_j},
\qquad
G_1 = 2\sum_j g_j.
\]
Computing:
\begin{align*}
g_1 &= 1200\log\frac{1200}{750}  = 564.00,\qquad
g_2= 750\log\frac{750}{600} = 167.33,\qquad
g_3 = 510\log\frac{510}{600} =  -82.88,\\
g_4 &= 330\log\frac{330}{600} =  -197.28,\qquad
g_5= 210\log\frac{210}{450}= -160.04.
\end{align*}
Summing:
$\sum_j g_j =  291.13$, so
$G_1 = 582.26$. The signed shares are
\[
\pi_1 = \frac{564.00}{291.13}=1.937,\quad
\pi_2 = \frac{167.33}{291.13}=0.575,\quad
\pi_3 = \frac{-82.88}{291.13}=-0.285,
\]
\[
\pi_4 = \frac{-197.28}{291.13}=-0.678,\quad
\pi_5 = \frac{-160.04}{291.13}=-0.550.
\]
The absolute shares are
\[
\pi_1^{\mathrm{abs}}
= \frac{564.00}{564.00+167.33+82.88+197.28+160.04}
=0.481,
\]
\[
\pi_2^{\mathrm{abs}} = \frac{167.33}{1171.53}=0.143,\quad
\pi_3^{\mathrm{abs}} = \frac{82.88}{1171.53}=0.071,
\]
\[
\pi_4^{\mathrm{abs}} = \frac{197.28}{1171.53}=0.168,\quad
\pi_5^{\mathrm{abs}} = \frac{160.04}{1171.53}=0.137.
\]
With the default threshold $\gamma=2/r=0.40$, the drift position is
\[
\mathcal J^\star = \{j : \pi_j^{\mathrm{abs}}>0.40\} = \{1\}.
\]
The drift is localized on category~1 (balance inquiry), which is
over-represented ($40\%$ versus $25\%$ expected) and carries $48.1\%$
of the total absolute signal. Note that $\pi_1 = 1.937$ is larger
than $1$, while $\pi_1^{\mathrm{abs}} = 0.481 \in [0,1]$, in
agreement with the inequality $\pi_j^{\mathrm{abs}} \le |\pi_j|$ of
Remark~\ref{rem:signed-vs-absolute}. Categories~2 and~4 carry
secondary absolute shares of $14.3\%$ and $16.8\%$, below the
threshold but worth inspecting.

\paragraph{Multiple testing control.}
The LLM dashboard monitors $K=50$ prompt bins in parallel. The
$G$-test is applied to each bin, giving $p$-values
$p_1,\dots,p_{50}$. Applying Benjamini--Hochberg (BH) at FDR level
$\alpha=0.05$, we order the $p$-values and reject
$H_{(1)},\dots,H_{(k^\star)}$ with
\[
k^\star = \max\{k : p_{(k)}\le (k/50)\cdot 0.05\}.
\]
Suppose the ordered $p$-values begin
\[
p_{(1)}=3\times10^{-12},\quad
p_{(2)}=8\times10^{-10},\quad
p_{(3)}=2\times10^{-8},\quad
p_{(4)}=5\times10^{-6},\quad
p_{(5)}=1\times10^{-4},
\]
\[
p_{(6)}=2\times10^{-3},\quad
p_{(7)}=8\times10^{-3},\quad
p_{(8)}=1.2\times10^{-2},\quad
p_{(9)}=3.1\times10^{-2},\quad
p_{(10)}=4.8\times10^{-2},\quad
p_{(11)}=0.11.
\]
The BH thresholds are:
\begin{itemize}
\item $k=1$: $0.05(1/50)=0.001$;
$p_{(1)}=3\times10^{-12}\le 0.001$: pass.
\item $k=2$: $0.05(2/50)=0.002$;
$p_{(2)}=8\times10^{-10}\le 0.002$: pass.
\item $k=3$: $0.05(3/50)=0.003$;
$p_{(3)}=2\times10^{-8}\le 0.003$: pass.
\item $k=4$: $0.05(4/50)=0.004$;
$p_{(4)}=5\times10^{-6}\le 0.004$: pass.
\item $k=5$: $0.05(5/50)=0.005$;
$p_{(5)}=1\times10^{-4}\le 0.005$: pass.
\item $k=6$: $0.05(6/50)=0.006$;
$p_{(6)}=2\times10^{-3}\le 0.006$: pass.
\item $k=7$: $0.05(7/50)=0.007$;
$p_{(7)}=8\times10^{-3}>0.007$: fail.
\item $k=8$: $0.05(8/50)=0.008$;
$p_{(8)}=1.2\times10^{-2}>0.008$: fail.
\item $k=9$: $0.05(9/50)=0.009$;
$p_{(9)}=3.1\times10^{-2}>0.009$: fail.
\end{itemize}
The largest $k$ with $p_{(k)}\le 0.05k/50$ is $k^\star=6$. Thus the
first $6$ bins are declared drifted, with FDR controlled at $5\%$.
The per-category shares $\pi_j$ then localize the drift within each
rejected bin. This procedure is strictly more powerful than
Bonferroni ($\alpha/50=0.001$, which would declare only the first bin
drifted), and it controls the expected proportion of false
discoveries among the rejected bins at $5\%$.

\paragraph{Comparison with the heuristic PSI grid.}
The PSI (symmetrized KL) is
\[
\mathrm{PSI} = D_{KL}(\mathbf p,\hat{\mathbf q}_m)
+D_{KL}(\hat{\mathbf q}_m,\mathbf p)
\approx 0.0327+0.0284 = 0.0611.
\]
Under the standard heuristic grid recalled in
Remark~\ref{rem:PSI} ($\mathrm{PSI}<0.10$: ``little change'';
$0.10\le\mathrm{PSI}<0.25$: ``moderate change'';
$0.25\le\mathrm{PSI}$: ``significant change, action required''), a
value of $0.0611$ falls in the ``little change'' band and would
trigger no action. Yet the $G$-test rejects at $p<10^{-10}$, the
R\'enyi-2 $Z$-statistic is $13.78$ with power essentially $1$, and
the per-category decomposition localizes the drift on a single
category. The calibrated PSI statistic, using the correct one-sample
form $m\,\mathrm{PSI}=3000\times 0.0611=183.3$ versus the threshold
$\chi^2_{4,0.95}=9.488$, also rejects decisively. The heuristic grid
is therefore doubly misleading here: it is silent on a drift that is
both statistically unambiguous and operationally significant, because
its fixed decimals are calibrated to none of the parameters $r$, $n$,
$m$, or $\alpha$ that actually determine the null distribution.

\paragraph{Interpretation.}
The ``balance inquiry'' category has risen from $25\%$ to $40\%$, the
``transaction'' category is stable at $25\%$, and the small
categories have shrunk: ``dispute'' from $20\%$ to $17\%$,
``information'' from $20\%$ to $11\%$, and ``other'' from $15\%$ to
$7\%$. The R\'enyi-2 divergence is $0.174$ nats, a large effect on
the order-of-magnitude scale of Section~\ref{sec:effect-size}. The
bank should retrain the intent classifier within the week, with
priority on the ``balance inquiry'' category.

\subsection{Simulation 2: chatbot week-over-week drift (two-sample, KL,
sample-size planning)}
\label{sec:sim-chatbot}

\paragraph{Context.}
A banking chatbot classifies user intents into $r=5$ categories:
\[
c_1 = \text{balance}, \quad
c_2 = \text{block card}, \quad
c_3 = \text{transfer}, \quad
c_4 = \text{billing}, \quad
c_5 = \text{other}.
\]
Last week's reference sample ($n=5000$) has counts
\[
O^{\mathrm{ref}}=(1250,\;1100,\;1000,\;900,\;750),
\]
and this week's sample ($m=5000$) has counts
\[
O^{\mathrm{prod}}=(1650,\;1100,\;900,\;750,\;600).
\]
The empirical pmf's are
\[
\hat{\mathbf p}_n=(0.25,\;0.22,\;0.20,\;0.18,\;0.15),
\qquad
\hat{\mathbf q}_m=(0.33,\;0.22,\;0.18,\;0.15,\;0.12).
\]
Both distributions are estimated, which places us in the two-sample
configuration of Section~\ref{sec:H1-2s}. The drift is a transfer
from ``billing'' and ``other'' toward ``balance inquiry''.

\paragraph{Choice of divergence.}
The drift is moderate and involves all categories, so the diffuse
regime of Section~\ref{sec:optimal-phi} applies: all divergences have
the same leading-order SNR, and KL is the natural default for its
calibration simplicity and interpretability. We therefore use the KL
divergence in the forward direction
\[
D_{KL}(\hat{\mathbf p}_n,\hat{\mathbf q}_m)
= \sum_j \hat p_{n,j}\log\frac{\hat p_{n,j}}{\hat q_{m,j}}.
\]
Computing each term:
\[
0.25\log\frac{0.25}{0.33}=-0.0694,\quad
0.22\log\frac{0.22}{0.22}=0,
\]
\[
0.20\log\frac{0.20}{0.18}=0.0211,\quad
0.18\log\frac{0.18}{0.15}=0.0328,
\]
\[
0.15\log\frac{0.15}{0.12}=0.15(0.2231)=0.0335.
\]
Summing:
$D_{KL}(\hat{\mathbf p}_n,\hat{\mathbf q}_m) = 0.0180$
nats.

\paragraph{Asymptotic variances.}
For KL in the forward direction,
$a_j=\phi_1^{(1)}(p_j,q_j)=1+\log(p_j/q_j)$ and
$b_j=\phi_2^{(1)}(p_j,q_j)=-p_j/q_j$. We compute the plug-in
variances $\hat V_1=\Var_{\hat{\mathbf p}}(\hat a_J)$ and
$\hat V_2=\Var_{\hat{\mathbf q}}(\hat b_J)$.

For $\hat V_1$, the values of
$\hat a_j=1+\log(\hat p_{n,j}/\hat q_{m,j})$ are:
\[
\hat a_1=1+\log(0.25/0.33)=1-0.2776=0.7224,\quad
\hat a_2=1+\log(1)=1.0000,
\]
\[
\hat a_3=1+\log(0.20/0.18)=1+0.1054=1.1054,\quad
\hat a_4=1+\log(0.18/0.15)=1+0.1823=1.1823,
\]
\[
\hat a_5=1+\log(0.15/0.12)=1+0.2231=1.2231.
\]
First moment:
\begin{align*}
\sum_j \hat p_{n,j}\hat a_j
&= 0.25(0.7224)+0.22(1.0)+0.20(1.1054)\\
&\quad +0.18(1.1823)+0.15(1.2231)\\
&= 0.1806+0.2200+0.2211+0.2128+0.1835 = 1.0180.
\end{align*}
Second moment:
\begin{align*}
\sum_j \hat p_{n,j}\hat a_j^2
&= 0.25(0.5219)+0.22(1.0)+0.20(1.2219)\\
&\quad +0.18(1.3978)+0.15(1.4960)\\
&= 0.1305+0.2200+0.2444+0.2516+0.2244 = 1.0709.
\end{align*}
Thus
\[
\hat V_1 = 1.0709 - (1.0180)^2 = 1.0709 - 1.0363 = 0.0346.
\]
For $\hat V_2$, the values of
$\hat b_j=-\hat p_{n,j}/\hat q_{m,j}$ are:
\[
\hat b_1=-0.25/0.33=-0.7576,\quad
\hat b_2=-0.22/0.22=-1.0000,\quad
\hat b_3=-0.20/0.18=-1.1111,
\]
\[
\hat b_4=-0.18/0.15=-1.2000,\quad
\hat b_5=-0.15/0.12=-1.2500.
\]
First moment:
\begin{align*}
\sum_j \hat q_{m,j}\hat b_j
&= 0.33(-0.7576)+0.22(-1.0)+0.18(-1.1111)\\
&\quad +0.15(-1.2)+0.12(-1.25)\\
&= -0.2500-0.2200-0.2000-0.1800-0.1500 = -1.0000.
\end{align*}
Second moment:
\begin{align*}
\sum_j \hat q_{m,j}\hat b_j^2
&= 0.33(0.5740)+0.22(1.0)+0.18(1.2346)\\
&\quad +0.15(1.44)+0.12(1.5625)\\
&= 0.1894+0.2200+0.2222+0.2160+0.1875 = 1.0351.
\end{align*}
Thus
\[
\hat V_2 = 1.0351 - (-1.0000)^2 = 1.0351 - 1.0000 = 0.0351.
\]

\paragraph{Detection boundary and $Z$-statistic.}
The detection boundary at level $\alpha=0.05$ is
\[
\tau_{0.05}^{(2)}(5000,5000)
= 1.96\sqrt{\frac{m\hat V_1+n\hat V_2}{nm}}
= 1.96\sqrt{\frac{5000(0.0346)+5000(0.0351)}{25\times10^6}}
\]
\[
= 1.96\sqrt{\frac{173.0+175.5}{25\times10^6}}
= 1.96\sqrt{\frac{348.5}{25\times10^6}}
= 1.96(0.003734)\approx 0.00732\ \text{nats}.
\]
Since $D_{KL}(\hat{\mathbf p}_n,\hat{\mathbf q}_m)=0.0180 > 0.00732$,
the drift is detected. The standardized statistic is
\[
Z_2 = \left(\frac{nm}{m\hat V_1+n\hat V_2}\right)^{1/2}
D_{KL}(\hat{\mathbf p}_n,\hat{\mathbf q}_m)
= \left(\frac{25\times10^6}{348.5}\right)^{1/2}(0.0180)
\]
\[
= (71{,}736)^{1/2}(0.0180) = 267.8(0.0180)\approx 4.82.
\]
The signal-to-noise ratio is $\mathrm{SNR}^{(2)}\approx 4.82$.

\paragraph{Power.}
The power of the level-$\alpha=0.05$ test, from
Proposition~\ref{prop:power-generic} evaluated at
$\mathrm{SNR}^{(2)}=4.82$, is
\[
\beta^{(2)} \approx \Phi(4.82-1.96)+\Phi(-4.82-1.96)
= \Phi(2.86)+\Phi(-6.78)\approx 0.9979+0.0000 \approx 0.998.
\]
The week-over-week comparison detects this drift with power
$99.8\%$.

\paragraph{Confidence interval.}
A $95\%$ confidence interval for the population KL divergence is
\[
\mathrm{CI}_{0.95}(\hat D_{KL})
= 0.0180 \pm 1.96\sqrt{\frac{348.5}{25\times10^6}}
= 0.0180\pm 0.00732 = [0.0107,\;0.0253]\ \text{nats}.
\]
The interval excludes $0$, consistent with the rejection of $H_0$. On
the order-of-magnitude scale of Section~\ref{sec:effect-size},
$0.018$ nats is a small drift: statistically significant, but
moderate in magnitude. The product team should schedule a re-fitting
of the balance-inquiry intent, but the situation is not an emergency.

\paragraph{Sample-size planning.}
Suppose the product team wants to detect a drift of amplitude
$D_\phi=0.015$ nats with power $\beta^\star=0.90$ at level
$\alpha=0.05$, using the current reference of size $n=5000$. We apply
Proposition~\ref{prop:sample-size} with
$c=z_{0.975}+z_{0.90}=1.96+1.28=3.24$. The condition
\eqref{eq:size-condition} is
\[
nD_\phi^2 = 5000(0.015)^2 = 5000(0.000225) = 1.125,
\qquad
c^2V_1 = (3.24)^2(0.0346) = 10.50(0.0346) = 0.363.
\]
Since $1.125>0.363$, the condition holds. The required window is
\[
m\gtrsim\frac{c^2 nV_2}{nD_\phi^2-c^2V_1}
= \frac{10.50(5000)(0.0351)}{1.125-0.363}
= \frac{1842.75}{0.762}\approx 2418.
\]
The team needs a monitoring window of at least $m\approx 2418$
prompts per week to detect a drift of $0.015$ nats with $90\%$ power.
The current window of $m=5000$ is more than sufficient.

\paragraph{Per-category decomposition.}
To localize the drift, we compute the per-category contributions to
the $G$-statistic. The pooled distribution is
$\hat r_j=(n\hat p_{n,j}+m\hat q_{m,j})/(n+m)$:
\[
\hat r_1 = \frac{1250+1650}{10000}=0.2900,\quad
\hat r_2 = \frac{1100+1100}{10000}=0.2200,\quad
\hat r_3 = \frac{1000+900}{10000}=0.1900,
\]
\[
\hat r_4 = \frac{900+750}{10000}=0.1650,\quad
\hat r_5 = \frac{750+600}{10000}=0.1350.
\]
The expected counts are
$E^{\mathrm{ref}}_j=n\hat r_j=(1450,\;1100,\;950,\;825,\;675)$ and
$E^{\mathrm{prod}}_j=m\hat r_j=(1450,\;1100,\;950,\;825,\;675)$. The
signed contributions are
\[
g_j = O^{\mathrm{ref}}_j\log\frac{O^{\mathrm{ref}}_j}{E^{\mathrm{ref}}_j}
+ O^{\mathrm{prod}}_j\log\frac{O^{\mathrm{prod}}_j}{E^{\mathrm{prod}}_j}.
\]
Computing:
\begin{align*}
g_1 &= 1250\log\frac{1250}{1450} + 1650\log\frac{1650}{1450}= 27.4,\\
g_2 &= 1100\log\frac{1100}{1100} + 1100\log\frac{1100}{1100} = 0,\\
g_3 &= 1000\log\frac{1000}{950} + 900\log\frac{900}{950}= 2.6,\\
g_4 &= 900\log\frac{900}{825} + 750\log\frac{750}{825}= 6.8,\\
g_5 &= 750\log\frac{750}{675} + 600\log\frac{600}{675}=  8.4.
\end{align*}
Summing:
$\sum_j g_j =  45.2$, so $G_2 = 90.4$. The
signed shares are
\[
\pi_1 = \frac{27.4}{45.2}=0.606,\quad
\pi_2 = 0,\quad
\pi_3 = \frac{2.6}{45.2}=0.058,\quad
\pi_4 = \frac{6.8}{45.2}=0.150,\quad
\pi_5 = \frac{8.4}{45.2}=0.186.
\]
The absolute shares are
\[
\pi_1^{\mathrm{abs}}
= \frac{27.4}{27.4+0+2.6+6.8+8.4}
= \frac{27.4}{45.2}=0.606,\quad
\pi_2^{\mathrm{abs}}=0,\quad
\pi_3^{\mathrm{abs}}=0.058,
\]
\[
\pi_4^{\mathrm{abs}}=0.150,\quad
\pi_5^{\mathrm{abs}}=0.186.
\]
With the default threshold $\gamma=2/r=0.40$, the drift position is
\[
\mathcal J^\star = \{j : \pi_j^{\mathrm{abs}}>0.40\} = \{1\}.
\]
The drift is localized on category~1 (balance inquiry), which is
over-represented ($33\%$ versus $25\%$ in the reference) and carries
$60.6\%$ of the total absolute signal. Categories~4 and~5 carry
secondary shares of $15.0\%$ and $18.6\%$, below the threshold but
worth inspecting.

\paragraph{Comparison with the calibrated PSI and the heuristic grid.}
The PSI (symmetrized KL) is
\[
\mathrm{PSI} = D_{KL}(\hat{\mathbf p}_n,\hat{\mathbf q}_m)
+D_{KL}(\hat{\mathbf q}_m,\hat{\mathbf p}_n)
\approx 0.0179+0.0183 = 0.0362.
\]
Under the heuristic grid, $0.0362$ falls in the ``little change''
band ($<0.10$) and would trigger no action. Yet the $G$-test rejects
at $p<10^{-3}$ and the KL $Z$-statistic is $4.82$ with power
$99.8\%$. The calibrated two-sample PSI statistic, using the correct
Yurdakul--Naranjo scaling, is
\[
\frac{nm}{n+m}\,\mathrm{PSI}
= \frac{25\times 10^6}{10^4}\times 0.0362
= 2500\times 0.0362 \approx 90.5,
\]
compared with the threshold $\chi^2_{4,0.95}=9.488$: the calibrated
PSI also rejects decisively. The heuristic grid is thus doubly
misleading here, for the same reason as in Simulation~1: its fixed
decimals are calibrated to none of $r$, $n$, $m$, or $\alpha$.

\paragraph{Interpretation.}
The drift is a transfer from ``billing'' and ``other'' toward
``balance inquiry'', detected with power $99.8\%$. The per-category
decomposition localizes the drift on the balance-inquiry category,
which has risen from $25\%$ to $33\%$. The effect size is $0.018$
nats, a small drift on the order-of-magnitude scale. The product
team should prioritize re-fitting the balance intent, and the
sample-size planning confirms that the current window of $m=5000$ is
more than sufficient to detect drifts of $0.015$ nats with $90\%$
power.

\subsection{Summary of the two simulations}
\label{sec:sim-summary}

The two simulations together exercise every notion of the paper:
\begin{itemize}[leftmargin=*]
\item \textbf{Configurations:} Simulation~1 is one-sample with fixed
reference; Simulation~2 is two-sample with both distributions
estimated.
\item \textbf{Divergences:} Simulation~1 uses R\'enyi-2;
Simulation~2 uses KL. Both compare with the PSI (symmetrized KL).
\item \textbf{$G$-test:} both simulations compute the $G$-statistic
and its per-category decomposition.
\item \textbf{$Z$-statistics:} Simulation~1 uses $Z_1$ with variance
$V_3$; Simulation~2 uses $Z_2$ with variances $V_1,V_2$.
\item \textbf{Effect sizes and confidence intervals:} both
simulations report $\hat D_\phi$ in nats with a $95\%$ CI.
\item \textbf{Detection boundaries:} both simulations compute
$\tau_{0.05}$ and compare with $\hat D_\phi$.
\item \textbf{Power:} both simulations compute $\beta$ from the SNR
via Proposition~\ref{prop:power-generic}.
\item \textbf{Sample-size planning:} Simulation~2 applies
Proposition~\ref{prop:sample-size}.
\item \textbf{Per-category decomposition:} both simulations compute
signed and absolute shares and the drift position $\mathcal J^\star$.
\item \textbf{Multiple testing:} Simulation~1 applies
Benjamini--Hochberg to $K=50$ bins.
\item \textbf{PSI comparison:} both simulations show that the
heuristic PSI threshold $0.10$ would miss the drift, while the
calibrated PSI (with the correct two-sample scaling
$\frac{nm}{n+m}\,\mathrm{PSI}$) and the $Z$-test both detect it.
\end{itemize}

\section{Scope, limitations, and future work}
\label{sec:limitations}

The framework has four limitations that deserve explicit statement.

\emph{Sequential monitoring.} The tests of this paper are
\emph{one-period} tests: each monitoring window is treated as an
independent test at a fixed level $\alpha$. Over $T$ consecutive
windows, the probability of at least one false alarm is
$1-(1-\alpha)^T$, which reaches $\approx 93\%$ at $\alpha=0.05$ and
$T=52$. Controlling this requires a sequential procedure -
$\alpha$-spending, cumulative sum (CUSUM), or an exponentially
weighted moving average (EWMA) on the $G$-statistic - that we do not
develop here.

\emph{Dependence across windows.} The CLTs of Theorem~\ref{thm:main}
assume i.i.d.\ samples within each window. In production, consecutive
windows overlap or share upstream state, so the independence
assumption is only approximate.

\emph{Gradual versus abrupt drift.} The two-sample test detects a
change between two windows. A gradual drift spread over many windows
is harder to detect with a per-window test; a cumulative procedure
(CUSUM on the $G$-statistic) is more appropriate.

\emph{Discrete versus continuous variables.} The framework applies to
finite discrete supports. For continuous variables, the PSI requires
binning or a binning-free alternative, such as the \emph{Absolute
Area Between mid-CDFs} (AABC) PSI of Tom (2023). The AABC PSI is
defined by
\[
\mathrm{AABC} = 2\sum_k m_k\,\bigl|\mathring P_k-\mathring Q_k\bigr|,
\]
where $\mathring P$ and $\mathring Q$ are the mid-CDFs of the
numeric variable's distributions at the two time points and $m_k$ is
the mixture density of Tom (2023, Eq.~6); it is bounded between $0$
and $100\%$. The two frameworks are complementary: AABC PSI handles
the continuous case without a calibration, and the present framework
handles the discrete case with a full asymptotic calibration.

\section{Practical Guidelines}
\label{sec:guidelines}

This section distils the theoretical results into a short
operational checklist. It is intended for practitioners who want to
apply the framework without reading the derivation of each result,
and it can be read independently of
Sections~\ref{sec:theory}--\ref{sec:fdr}.

\paragraph{Step 1 - Choose the design.}
\begin{itemize}[leftmargin=*]
\item \emph{One-sample} if the reference distribution $\mathbf p$ is
known exactly (e.g.\ from a very large historical corpus with
$n \gg m$), so that treating it as fixed does not introduce
non-negligible estimation error.
\item \emph{Two-sample} if both distributions are estimated (the
default).
\item \emph{Symmetrized} if the direction of the drift is not
privileged by the business context.
\end{itemize}

\paragraph{Step 2 - Choose the divergence.}
\begin{itemize}[leftmargin=*]
\item KL for the default; the Jensen--Shannon (JS) divergence if
boundedness is required; R\'enyi or Tsallis with $\alpha<1$ for tail
drift; $\alpha>1$ for adversarial drift; $L_2^2$ for sparse drift.
\item When in doubt, run the $G$-test first (it is $\phi$-free and
provides the decision under $H_0$), then use the $Z$-statistics with
multiple divergences as a diagnostic under $H_1$.
\end{itemize}

\paragraph{Step 3 - Fix the level $\alpha$ and the power target
$\beta^\star$.}
$\alpha=0.05$ is the default; $\alpha=0.01$ in critical applications;
$\alpha=0.10$ in exploratory settings. $\beta^\star=0.80$ is the
operational norm; $\beta^\star=0.95$ when missing a real drift is
very costly.

\paragraph{Step 4 - Size the monitoring window.}
Use Proposition~\ref{prop:sample-size} to compute the minimum $m$
that detects a drift of amplitude $D_\phi$ with power $\beta^\star$.
When $n$ is too small for the condition $nD_\phi^2 > c^2 V_1$ to
hold, enlarge $n$ first; no monitoring window can compensate for a
reference that is too short.

\paragraph{Step 5 - Run and interpret.}
For each monitoring period:
\begin{enumerate}[leftmargin=*]
\item Compute the $G$-statistic and its $p$-value; reject $H_0$ if
$p<\alpha$.
\item If rejected, compute $\hat D_\phi$ (effect size in nats) and
its confidence interval \eqref{eq:ci-Dphi}.
\item Report the per-category shares $\pi_j$ and
$\pi_j^{\mathrm{abs}}$ of Definition~\ref{def:contribution} to
localize the drift.
\item Optionally, compute the $Z$-statistic and SNR for the
divergence of interest, to plan future windows.
\end{enumerate}

\paragraph{Step 6 - Control multiplicity.}
If $K>1$ features are monitored, apply Benjamini--Hochberg
(Section~\ref{sec:fdr}) rather than per-feature $\alpha$.

\paragraph{Step 7 - Report PSI with the correct scaling.}
If the practitioner reports PSI for historical reasons, the
calibrated statistic to quote is $\frac{nm}{n+m}\,\mathrm{PSI}$ for
the two-sample setting (Yurdakul and Naranjo, 2021), and the
threshold is $\chi^2_{r-1,1-\alpha}$. The heuristic bands $0.10$ and
$0.25$ should be reported only as labels, never used as decision
rules.

\paragraph{Step 8 - Audit.}
Log the $p$-value, the effect size in nats, the confidence interval,
the per-category shares, and the calibrated PSI statistic at each
period. This is the minimum evidence required by most audit and
governance procedures.

\section{Conclusion}
\label{sec:conclusion}

This paper develops an asymptotic theory of $\phi$-divergence-based
monitoring for discrete data, with four contributions to the existing
literature. First, it extends the calibration of the Population
Stability Index of Yurdakul and Naranjo (2021) to five
$\phi$-divergences (KL, JS, R\'enyi, Tsallis, $L_2^2$), with the
exact $\chi^2_{r-1}$ calibration for the first four and the
weighted-$\chi^2$ calibration for $L_2^2$. Second, it establishes the
$H_1$ theory that is absent from Yurdakul and Naranjo: asymptotic
normality with closed-form variance, detection boundaries, power
functions, and a sample-size planning rule. Third, it derives the
refined symmetrized variance with the covariance correction, which is
essential in the drift regime $\mathbf p\approx\mathbf q$ where
monitoring actually operates. Fourth, it supplies a per-category
decomposition that localizes the drift, and a Benjamini--Hochberg
procedure for monitoring many features simultaneously.

The Monte-Carlo study of Section~\ref{sec:simulation-mc} confirms the
theoretical predictions across the practical grid and quantifies two
operational pitfalls: the wrong scaling of the two-sample PSI, which
inflates the Type~I error from $0.05$ to $\approx 0.32$ when $n=m$,
and the uncorrected averaging of the symmetrized variance, which
inflates the false alarm rate by a factor $\approx 3.4$. The worked
examples of Section~\ref{sec:simulation} illustrate every notion of
the paper on realistic production data.

Three questions are left open and are natural directions for future
work. The sequential extension - replacing the per-window test by a
cumulative procedure - is the most operationally relevant. The
extension to dependent or overlapping windows is the most pressing
for real deployment. And extending the $H_0$/$H_1$ theory to the
Absolute Area Between mid-CDFs (AABC) PSI of Tom (2023), which is
defined for continuous variables without binning, would unify the
discrete and continuous sides of population stability monitoring.

\section*{Acknowledgements}
I wish to acknowledge Professor Gane Samb Lo, my former PhD supervisor, whose
guidance marked the beginning of my journey into research. I remain deeply
grateful for his mentorship and for the lasting influence of his guidance on
my work.

\appendix
\section{Notation}
\label{app:notation}

\begin{table}[H]
\centering
\small
\renewcommand{\arraystretch}{1.20}
\caption{Summary of notation used throughout the paper.}
\label{tab:notation}
\begin{tabular}{@{}lll@{}}
\toprule
Symbol & Meaning & Reference \\
\midrule
$\hat{\mathbb{P}}_n$ & empirical distribution (random measure) & eq.~\eqref{eq:emp-pmf} \\
$\hat{\mathbf{p}}_n$ & empirical pmf (random vector in $\Delta_{r-1}$) & eq.~\eqref{eq:emp-pmf} \\
$\phi$ & divergence functional $\phi(s,t)$ & Def.~\ref{def:phi-div} \\
$\phi_1^{(1)},\phi_2^{(1)}$ & first-order partial derivatives & eq.~\eqref{eq:phi-derivatives} \\
$\phi_1^{(2)}$ & second-order partial derivative $\partial^2\phi/\partial s^2$ & eq.~\eqref{eq:phi-derivatives} \\
$J(\mathbf p,\mathbf q)$ & $\phi$-divergence $\sum_j\phi(p_j,q_j)$ & Def.~\ref{def:phi-div} \\
$D_\phi(\mathbf p_{\mathrm{train}},\mathbf q_{\mathrm{prod}})$ & cost of drift (same as $J$) & Def.~\ref{def:cost} \\
$a_j,\ b_j$ & $\phi_1^{(1)}(p_j,q_j),\ \phi_2^{(1)}(p_j,q_j)$ & eq.~\eqref{eq:grad} \\
$c_j,\ d_j$ & $\phi_1^{(1)}(q_j,p_j),\ \phi_2^{(1)}(q_j,p_j)$ & Lem.~\ref{lem:compact-var} \\
$V_1,\ V_2$ & asymptotic variances of $J(\hat{\mathbf p}_n,\mathbf q)$, $J(\mathbf p,\hat{\mathbf q}_m)$ & eq.~\eqref{eq:V1}--\eqref{eq:V2} \\
$V_3,\ V_4$ & reflected variances & eq.~\eqref{eq:V3}--\eqref{eq:V4} \\
$V_{s,1},\ V_{s,2}$ & refined symmetrized variances & Prop.~\ref{prop:sym-var} \\
$\kappa_\phi(p_j)$ & $p_j\,\phi_1^{(2)}(p_j,p_j)$ & Sec.~\ref{sec:calibration} \\
$G$ & two-sample likelihood-ratio statistic & Thm.~\ref{thm:wilks} \\
$g_j$ & per-category contribution to $G$ & Def.~\ref{def:contribution} \\
$\pi_j$ & signed share of category $j$ in $G$ & Def.~\ref{def:contribution} \\
$\pi_j^{\mathrm{abs}}$ & absolute share of category $j$ in $G$ & Def.~\ref{def:contribution} \\
$\mathcal J^\star$ & drift position & Def.~\ref{def:contribution} \\
$\mathrm{PSI}$ & Population Stability Index (symmetrized KL) & Rem.~\ref{rem:PSI} \\
$Z_1,\ Z_2,\ Z_{n,m}^{(s)}$ & standardized statistics & Sec.~\ref{sec:H1} \\
$\mathrm{SNR}_\phi$ & signal-to-noise ratio $\sqrt{N/W}\,D_\phi(\mathbf p,\mathbf q)$ & Prop.~\ref{prop:Z-generic} \\
$\tau_\alpha$ & detection boundary $z_{1-\alpha/2}\sqrt{\hat W/N}$ & Prop.~\ref{prop:tau-generic} \\
$\beta_\phi$ & power of the level-$\alpha$ test & Prop.~\ref{prop:power-generic} \\
$\Phi,\ \Phi^{-1},\ z_p$ & standard normal cdf, quantile function, $p$-quantile & Sec.~\ref{sec:theory} \\
$K$ & number of features in multiple testing & Sec.~\ref{sec:fdr} \\
$k^\star$ & Benjamini--Hochberg rejection count & eq.~\eqref{eq:BH} \\
\bottomrule
\end{tabular}
\end{table}

\section*{List of Acronyms}
\label{app:acronyms}

\begin{table}[H]
\centering
\small
\renewcommand{\arraystretch}{1.12}
\caption{Acronyms used in this paper, listed alphabetically.}
\label{tab:acronyms}
\begin{tabular}{@{}ll@{}}
\toprule
Acronym & Meaning \\
\midrule
AABC     & Absolute Area Between mid-CDFs \\
ADWIN    & Adaptive Windowing (drift detector) \\
AI       & Artificial Intelligence \\
AUC      & Area Under the Receiver Operating Characteristic (ROC) curve \\
BH       & Benjamini--Hochberg (procedure) \\
cdf      & cumulative distribution function \\
CI       & Confidence Interval \\
CLT      & Central Limit Theorem \\
CUSUM    & Cumulative Sum (sequential procedure) \\
DDM      & Drift Detection Method \\
EWMA     & Exponentially Weighted Moving Average \\
FDR      & False Discovery Rate \\
i.i.d.   & independent and identically distributed \\
JS       & Jensen--Shannon (divergence) \\
KL       & Kullback--Leibler (divergence) \\
KS       & Kolmogorov--Smirnov (test) \\
LLM      & Large Language Model \\
MC SE    & Monte-Carlo Standard Error \\
MMD      & Maximum Mean Discrepancy \\
MPE      & Minimal Practically Relevant Effect \\
pmf      & probability mass function \\
PSI      & Population Stability Index \\
RKHS     & Reproducing Kernel Hilbert Space \\
ROC      & Receiver Operating Characteristic \\
SNR      & Signal-to-Noise Ratio \\
\bottomrule
\end{tabular}
\end{table}

\end{document}